\documentclass[review]{elsarticle}
\usepackage{lineno}
\usepackage{hyperref}
\usepackage{framed,multirow}
\usepackage{verbatim}
\usepackage{amssymb}
\usepackage{latexsym}
\usepackage{amsmath}
\usepackage{lipsum}
\usepackage{amsfonts}
\usepackage{graphicx}
\usepackage{epstopdf}
\usepackage{bm,nicefrac}
\usepackage{subfig}
\usepackage{algorithm,algpseudocode,float}
\usepackage{setspace}
\usepackage{algcompatible}
\usepackage{tikz}
\usepackage{pgfplots}
\usepackage{mathptmx}
\usepackage{amsmath}
\usepackage{amssymb}
\usepackage{amsthm}
\usepackage{amsmath}
\usepackage{amsthm}
\usepackage{amssymb}
\usepackage{listings}
\usepackage{color}
\usepackage[justification=centering]{caption}
\usepackage[numbers]{natbib}
\usepackage{picture}
\usepackage{subfig}
\usepackage{float}
\usepackage{mathrsfs}
\usepackage{bm}
\usepackage{tikz}
\usepackage{pgfplots}
\usepackage{setspace}
\usepackage{algcompatible}
\usepackage{mathtools}
\usepackage{epstopdf}
\usepackage{tikz}
\usepackage{mathptmx}
\usepackage{stmaryrd}
\usepackage{algorithm,algpseudocode,float}
\usepackage{filecontents}
\usepackage{comment}

\newtheorem{thm}{Theorem}[section]

\theoremstyle{definition}

\newcommand{\p}{\partial}
\newcommand{\wh}{\widehat}
\newcommand{\jump}[1]{\ensuremath{\llbracket#1\rrbracket}}
\newcommand{\avg}[1]{\ensuremath{\LRc{\!\LRc{#1}\!}}}
\newcommand{\LRc}[1]{\left\{ #1 \right\}}

\usepackage{url}
\usepackage{xcolor}
\usepackage{comment}
\renewcommand{\hat}{\widehat}
\renewcommand{\tilde}{\widetilde}
\definecolor{newcolor}{rgb}{.8,.349,.1}
\newcommand{\logLogSlopeTriangle}[5]
{
	\pgfplotsextra
	{
		\pgfkeysgetvalue{/pgfplots/xmin}{\xmin}
		\pgfkeysgetvalue{/pgfplots/xmax}{\xmax}
		\pgfkeysgetvalue{/pgfplots/ymin}{\ymin}
		\pgfkeysgetvalue{/pgfplots/ymax}{\ymax}
		
		\pgfmathsetmacro{\xArel}{#1}
		\pgfmathsetmacro{\yArel}{#3}
		\pgfmathsetmacro{\xBrel}{#1-#2}
		\pgfmathsetmacro{\yBrel}{\yArel}
		\pgfmathsetmacro{\xCrel}{\xArel}
		\pgfmathsetmacro{\lnxB}{\xmin*(1-(#1-#2))+\xmax*(#1-#2)} % in [xmin,xmax].
		\pgfmathsetmacro{\lnxA}{\xmin*(1-#1)+\xmax*#1} % in [xmin,xmax].
		\pgfmathsetmacro{\lnyA}{\ymin*(1-#3)+\ymax*#3} % in [ymin,ymax].
		\pgfmathsetmacro{\lnyC}{\lnyA+#4*(\lnxA-\lnxB)}
		\pgfmathsetmacro{\yCrel}{\lnyC-\ymin)/(\ymax-\ymin)} % THE IMPROVED EXPRESSION WITHOUT 'DIMENSION TOO LARGE' ERROR.
		
		\coordinate (A) at (rel axis cs:\xArel,\yArel);
		\coordinate (B) at (rel axis cs:\xBrel,\yBrel);
		\coordinate (C) at (rel axis cs:\xCrel,\yCrel);
		
		\draw[#5]   (A)-- node[pos=0.5,anchor=north] {}%\tiny{1}}
		(B)-- 
		(C)-- node[pos=0.5,anchor=west] {\tiny{#4}}
		cycle;
	}
}

\journal{International Journal for Numerical Methods in Engineering}

\begin{document}

%\verso{Guo, Chan}

\begin{frontmatter}

\title{High order weight-adjusted discontinuous Galerkin methods for wave propagation on moving curved meshes}%
%\tnotetext[tnote1]{This is an example for title footnote coding.}

%\author[1]{Kaihang {Guo}\corref{cor1}}
%\author[1]{Jesse {Chan}}
%\cortext[cor1]{Corresponding author.\\
%Email addresses: Kaihang.Guo@rice.edu (K. Guo), Jesse.Chan@rice.edu (J. Chan).}
%\footnote{Email}{}
%\fntext[]{This is author footnote for second author.}  
%\author[2]{Given-name3 \snm{Surname3}}
%% Third author's email
%\ead{author3@author.com}
%\author[2]{Given-name4 \snm{Surname4}}
\author[1]{Kaihang Guo\corref{cor1}}
\ead{Kaihang.Guo@rice.edu}
\cortext[cor1]{Principal Corresponding author}
\author[1]{Jesse Chan}
\ead{Jesse.Chan@rice.edu}

\address[1]{Department of Computational and Applied Mathematics, Rice University, 6100 Main St, Houston, TX 77005, United States}
%\address[2]{Department of Pediatrics-Cardiology, Baylor College of Medicine, Houston, TX, United States}

%\received{1 May 2013}
%\finalform{10 May 2013}
%\accepted{13 May 2013}
%\availableonline{15 May 2013}
%\communicated{S. Sarkar}
\begin{abstract}
 This paper presents high order accurate discontinuous Galerkin (DG) methods for wave problems on moving curved meshes with general choices of basis and quadrature. The proposed method adopts an arbitrary Lagrangian-Eulerian (ALE) formulation to map the acoustic wave equation from the time-dependent moving physical domain onto a fixed reference domain. For moving curved meshes, weighted mass matrices must be assembled and inverted at each time step when using explicit time stepping methods. We avoid this step by utilizing an easily invertible weight-adjusted approximation. The resulting semi-discrete weight-adjusted DG scheme is provably energy stable up to a term which converges to zero with the same rate as the optimal $L^2$ error estimate. Numerical experiments using both polynomial and B-spline bases verify the high order accuracy and energy stability of proposed methods. 
\end{abstract}


%\begin{keyword}
% MSC codes here, in the form: \MSC code \sep code
% or \MSC[2008] code \sep code (2000 is the default)
%\MSC 41A05\sep 41A10\sep 65D05\sep 65D17
% Keywords
%\KWD discontinuous Galerkin\sep Bernstein\sep high order\sep heterogeneous media
%\end{keyword}
\end{frontmatter}

%\linenumbers
 
\section{Introduction}
Efficient and accurate simulations of wave propagation have a wide range of applications in science and engineering, from seismic and medical imaging to rupture and earthquake simulations. Moving meshes appear when simulating problems with moving domains or boundaries \cite{ostashev2015acoustics,morse1986theoretical}, e.g., wave scattering from moving or vibrating boundaries \cite{winters2014ale,winters2014discontinuous}. Another application of moving mesh methods is to resolve sharp wave fronts or localized features by clustering mesh grid points in certain area with large solution variations (dynamic $r$-adaptivity) \cite{budd2009adaptivity,cao2002moving,tang2005moving}. 

In this work, we use an arbitrary Lagrangian-Eulerian (ALE) formulation of partial differential equations (PDEs) to account for mesh motion. The ALE method \cite{noh1963cel,franck1964mixed,trulio1966theory,donea1977lagrangian} is a generalization of the Lagrangian and Eulerian methods, which allows mesh grid point move arbitrarily. The basic idea of the ALE transformation is to map governing equations from a time-dependent physical domain $\Omega_t$ to a fixed reference domain $\wh{\Omega}$. 

Recently, discontinuous Galerkin (DG) methods based on the ALE framework have been developed for their advantages of high order approximation on arbitrary unstructured meshes \cite{kopriva2016provably,wang2015high,fu2018arbitrary}. A variety of DG spectral element methods (DG-SEM) for wave propagation on moving domains can be found in \cite{kopriva2016provably,minoli2011discontinuous,winters2014discontinuous,winters2014ale}. However, it is known that DG methods for problems on curvilinear meshes \cite{kopriva2016geometry} can be unstable even with static domains, and care must be taken to ensure stability for moving curved meshes. Several approaches have been proposed to ensure energy stability of high order DG-type formulations. For constant coefficient hyperbolic systems, Nikkar and Nordstr$\ddot{\textmd{o}}$m \cite{nikkar2015fully} develop a high order, fully discrete, conservative and energy stable finite difference scheme using summation-by-parts (SBP) operators both in time and space. In \cite{kopriva2016provably}, Kopriva et al.\ proposed a provably stable semi-discrete DG-SEM approximation on moving hexahedral meshes by averaging the conservative form and a non-conservative form of the equation. However, this method requires the use of tensor product elements and tensor product Legendre-Gauss-Lobatto quadrature. For moving meshes with non-curved triangular elements, Fu et al.\ \cite{fu2018arbitrary} construct an ALE-DG method by assuming mesh motion is restricted to affine time-dependent mappings.

In this paper, we are interested in constructing more general stable DG formulations on moving curved meshes without restrictions on element type, quadrature, or choice of local approximation space. During mesh motion, affine elements may become curved elements, so the determinant of the Jacobian matrix given by the ALE transformation is time-dependent and spatially varying. In this case, we must construct and invert a high order weighted mass matrix on each element at every time step. In order to reduce the computational cost, we build upon a weight-adjusted DG (WADG) formulation \cite{chan2017weight}, which is low storage, energy stable, and high order accurate  for static heterogeneous media \cite{chan2017weight,guo2020bernstein} and curvilinear meshes \cite{chan2017curved,guo2020weight}. We then extend this WADG formulation to moving curved meshes, and prove that it is energy stable up to a term which converges to zero with the same rate as the optimal $L^2$ error estimate. We also introduce a novel penalty flux (a blend of upwind and Lax-Friedrichs fluxes \cite{Hesthaven2007}) for moving meshes. The proposed ALE-DG method can be applied on unstructured triangular and quadrilateral meshes using high order polynomial and non-polynomial bases. We implement numerical experiments using both polynomial and B-spline bases to verify the high order accuracy and energy stability of proposed methods. 

The reminder of this paper is organized as follows: In Section~\ref{sec:ale}, we derive an ALE formulation of a conservation law for the wave equation on moving meshes. In Section~\ref{sec:mesh-ale}, we use the constant solution on moving meshes to show results of energy conservation for the proposed ALE-DG methods using standard Galerkin approach and weight-adjusted approach, respectively. In Section~\ref{sec:aledg}, we propose an ALE-DG method for the acoustic wave equation with a novel penalty flux. In Section~\ref{sec:ale-experiment}, we present numerical experiments to verify energy stability and high order accurate of the proposed methods.

\section{Mathematical notation}\label{sec:notation}
In this section, we introduce mathematical notation that will be used in the following sections. In this paper, we restrict ourselves to two-dimensional problems, but all results can be directly extended to three dimensions. 

We assume that the fixed non-curved reference domain $\wh{\Omega}$ is exactly represented by a triangularization $\wh{\Omega}_h=\cup^{K}_{k=1}D^k$, which consists of $K$ non-overlapping elements. Each element $D^k$ is the image of the reference element $\hat{D}$ under an affine mapping $\bm{\Phi}^k$
$$\wh{\bm{x}}=\bm{\Phi}^k\bm{r},\qquad \wh{\bm{x}}\in D^k,\quad \bm{r}\in \hat{D},$$ 
where $\wh{\bm{x}}=\left(\xi_1,\xi_2\right)$ are coordinates on the $k$th element and $\bm{r}=\left(r,s\right)$ are coordinates on the reference element. In this work, the reference element  is taken to be bi-unit right triangle,
$$\wh{D}=\{\left(r,s\right)\geq-1,\quad r+s\leq 0\}.$$
Over each element $D^k$, we define the approximation space $V_h\left(D^k\right)$ as
$$V_h\left(D^k\right)= V_h\left(\hat{D}\right)\circ\left(\bm{\Phi}^k\right)^{-1}=\Big\{\widehat{v}_h\circ\left(\bm{\Phi}^k\right)^{-1},\  \widehat{v}_h\in V_h(\widehat{D}) \Big\},$$
where $V_h\left(\hat{D}\right)$ is a  polynomial approximation space of degree $N$ on the reference element defined by 
$$V_h\left(\widehat{D}\right)=P^N\left(\widehat{D}\right)=\big\{r^i s^j,\quad 0\leq i+j\leq N\big\}.$$
We will also consider quadrilateral meshes and non-polynomial approximation spaces such as splines in Section~\ref{sec:numerical-spline}.

For simplicity, we denote the $L^2$ inner product over $D^k$ and over surface of $\partial D^k$ as
\begin{align*}
\left(u,v\right)_{L^2(D^k)} = \int_{D^k}uv,\qquad \left\langle u,v\right\rangle_{L^2(\p D^k)} = \int_{\p  D^k}uv.
\end{align*}
We define the jump and average of a scalar variable $p$  across element interfaces as
$$\jump{p}=p^+-p, \qquad \avg{p}=\frac{1}{2}\left(p^++p\right),$$
where $p^+$ and $p$ are the neighboring and local traces of the solution over each face, respectively. The jump and average of a vector variable $\bm{u}$  are defined component-wise as follows
$$\left(\jump{\bm{u}}\right)_i=\jump{\bm{u}_i},\qquad \left(\avg{\bm{u}}\right)_i=\frac{1}{2}\avg{\bm{u}_i}.$$

\section{Arbitrary Lagrangian-Eulerian formulations} 
\label{sec:ale}
To take into account mesh motion and deformation, we use an arbitrary Lagrangian-Eulerian (ALE) formulation in this work. In this section, we derive an ALE formulation of conservation laws for the first-order system of two-dimensional acoustic wave equations. We assume constant wavespeed $c=1$ and write the acoustic wave equation into a conservative matrix form
\begin{equation}\label{eq:conservative}
\frac{d\bm{q}}{dt}+\sum_i\frac{\p\bm{f}_i}{\p x_i}=0,
\end{equation}
where
$$\bm{q}=\renewcommand\arraystretch{0.7}\begin{pmatrix}
p\\u\\v
\end{pmatrix}
,\quad \bm{f}_1=\begin{pmatrix}
0&1&0\\
1&0&0\\
0&0&0
\end{pmatrix}\begin{pmatrix}
p\\u\\v
\end{pmatrix},\quad
\bm{f}_2=\begin{pmatrix}
0&0&1\\
0&0&0\\
1&0&0
\end{pmatrix}\begin{pmatrix}
p\\u\\v
\end{pmatrix}.$$
\begin{figure}
	\centering
	\includegraphics[width=0.7\linewidth]{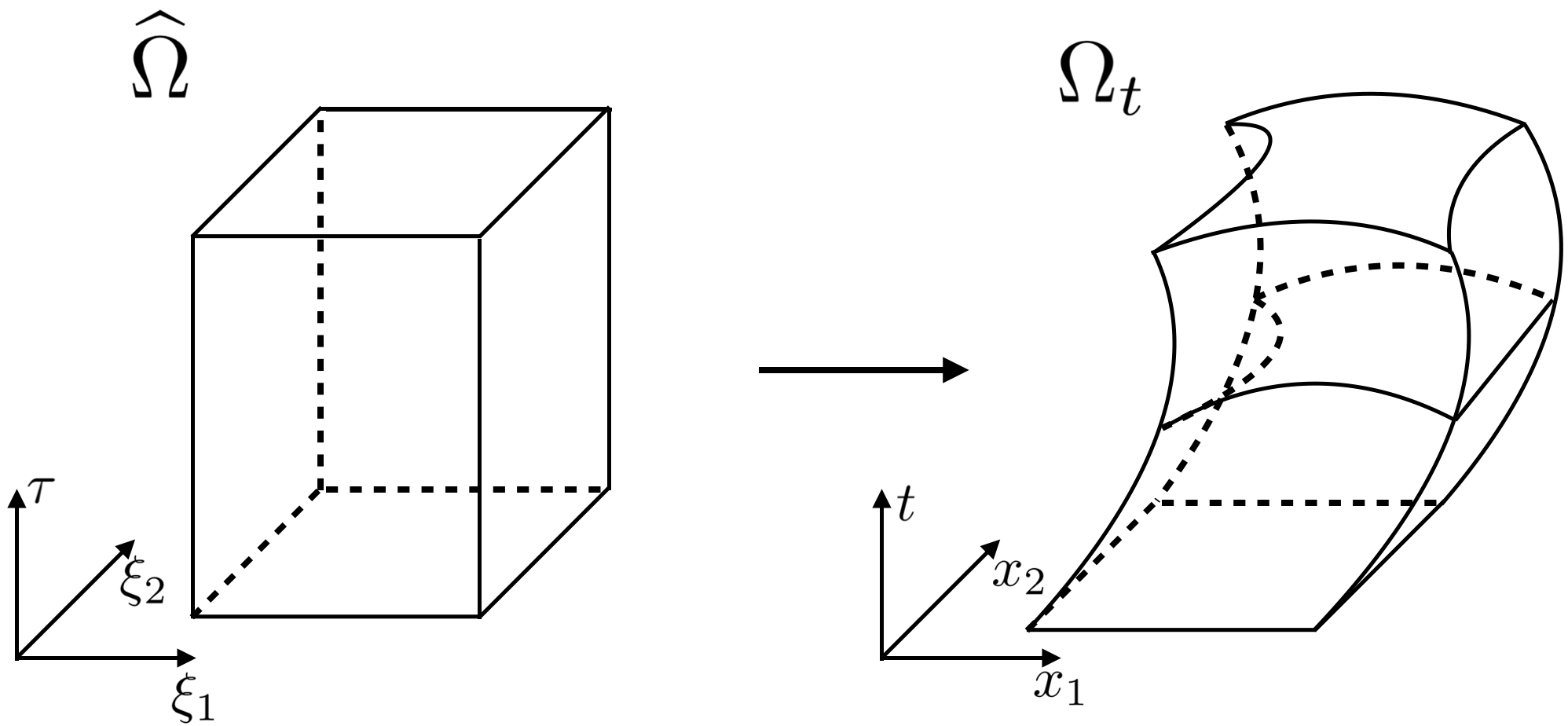}
	\caption{Sketch of the ALE transformation.}
	\label{fig:ALE}
\end{figure}

In the ALE formulation, we need to map the fixed reference domain $\wh{\bm{x}}=\left(\xi_1,\xi_2\right)\in \widehat{\Omega}$ with time $\tau$ onto the time-dependent physical domain $\bm{x}=\left(x_1,x_2\right)\in \Omega_t$ with time $t$. We assume this two coordinate systems are connected through a transformation as illustrated in Figure~\ref{fig:ALE}
$$t=\tau,\qquad \bm{x} = \bm{\textit{X}}\left(\wh{\bm{x}},\tau\right).$$ 
This transformation can be described by the following Jacobian matrix
$$\renewcommand\arraystretch{0.9}\frac{\partial\left(\bm{x},t\right)}{\partial\left(\bm{\xi},\tau\right)}=\begin{pmatrix}
\frac{\partial x_1}{\partial \xi_1} &\frac{\partial x_1}{\partial \xi_2}&\frac{\partial x_1}{\partial \tau}\\
\frac{\partial x_2}{\partial \xi_1} &\frac{\partial x_2}{\partial \xi_2}&\frac{\partial x_2}{\partial \tau}\\
0&0&1\\
\end{pmatrix}, \qquad J = \Bigg|\frac{\partial\left(\bm{x},t\right)}{\partial\left(\wh{\bm{x}},\tau\right)}\Bigg|.$$
From this transformation, we obtain the following relationship between differential operators on the physical element and on the reference domain 
\begin{equation}
\begin{split}
\frac{\partial}{\partial t}&=\frac{\partial}{\partial \tau}+\sum_{j}\frac{\partial\xi_j}{\partial t}\frac{\partial}{\partial\xi_j}\\
\frac{\partial}{\partial x_i}&=\sum_{j}\frac{\partial \xi_j}{\partial x_i}\frac{\partial}{\partial \xi_j},\qquad i =1,2.
\end{split}\label{eq:ale-relation}
\end{equation}
Now, we replace differential operators in (\ref{eq:conservative}) by (\ref{eq:ale-relation})
\begin{align}
\frac{d\bm{q}}{d\tau}+\sum_{j}\frac{\p \xi_j}{\p t}\frac{\p\bm{q}}{\p\xi_j}+\sum_i\sum_j\frac{\p\xi_j}{\p x_i}\frac{\p \bm{f}_i}{\p\xi_j}=0.
\end{align}
Multiplying by $J$, we obtain that
\begin{align}
J\frac{d\bm{q}}{d\tau}+J\sum_{j}\frac{\p \xi_j}{\p t}\frac{\p\bm{q}}{\p\xi_j}+J\sum_i\sum_j\frac{\p\xi_j}{\p x_i}\frac{\p \bm{f}_i}{\p\xi_j}=0.
\label{eq:ale1}
\end{align}
Then, we rewrite (\ref{eq:ale1}) as
\begin{align}
&\frac{d\bm{q}J}{d\tau}-\frac{dJ}{d\tau}\bm{q}+\sum_{j}\frac{\p \xi_j}{\p t}\frac{\p \bm{q}J}{\p\xi_j}-\sum_{j}\frac{\p \xi_j}{\p t}\frac{\p J}{\p\xi_j}\bm{q}+\sum_i\sum_j\frac{\p\xi_j}{\p x_i}\frac{\p J\bm{f}_i}{\p\xi_j}-\sum_i\sum_j\frac{\p\xi_j}{\p x_i}\frac{\p J}{\p\xi_j}\bm{f}_i\nonumber
\\
=&\frac{d\bm{q}J}{d\tau}+\sum_{j}\frac{\p \xi_j}{\p t}\frac{\p \bm{q}J}{\p\xi_j}+\sum_i\sum_j\frac{\p\xi_j}{\p x_i}\frac{\p J\bm{f}_i}{\p\xi_j}-\left(\frac{dJ}{d\tau}+\sum_{j}\frac{\p \xi_j}{\p t}\frac{\p J}{\p\xi_j}\right)\bm{q}-\sum_i\left(\sum_j\frac{\p\xi_j}{\p x_i}\frac{\p J}{\p\xi_j}\right)\bm{f}_i=0.
\label{eq:ale3}
\end{align}
For continuous geometric mappings, $J$ satisfies both the well-known metric identities \cite{mavriplis2006construction}
\begin{align}
\sum_j\frac{\p\xi_j}{\p x_i}\frac{\p J}{\p\xi_j}=0,
\label{eq:metric-id}
\end{align}
and the Geometric Conservation Law (GCL) \cite{vinokur1974conservation}
\begin{align}
\frac{dJ}{d\tau}+\sum_{j}\frac{\p \xi_j}{\p t}\frac{\p J}{\p\xi_j}=0.
\label{eq:gcl}
\end{align}
Simplifying (\ref{eq:ale3}) using (\ref{eq:metric-id}) and (\ref{eq:gcl}), we obtain an equivalent conservation law to (\ref{eq:conservative}) on $\widehat{\Omega}$
\begin{align}
\frac{d\bm{q}J}{d\tau}+\sum_{j}\frac{\p \xi_j}{\p t}\frac{\p \bm{q}J}{\p\xi_j}+\sum_i\sum_j\frac{\p\xi_j}{\p x_i}\frac{\p J\bm{f}_i}{\p\xi_j}=0.
\label{eq:ale4}
\end{align}
In two dimensions, (\ref{eq:ale4}) is equivalent to
\begin{equation}
\begin{split}
\renewcommand\arraystretch{0.7}\begin{pmatrix}
pJ\\
uJ\\
vJ
\end{pmatrix}_\tau + \left(\bm{A}_1\begin{pmatrix}
p\\u\\v
\end{pmatrix}\right)_{\xi_1}+\left(\bm{A}_2\begin{pmatrix}
p\\u\\v
\end{pmatrix}\right)_{\xi_2}=0,
\end{split}
\label{eq:ale-matrix}
\end{equation}
where
$$\bm{A}_1=\begin{pmatrix}
%\begin{smallmatrix}
\frac{\p\xi_1}{\p t}J & \frac{\p\xi_1}{\p x_1}J & \frac{\p\xi_1}{\p x_2}J\\
\frac{\p\xi_1}{\p x_1}J&\frac{\p\xi_1}{\p t}J&0\\
\frac{\p\xi_1}{\p x_2}J&0&\frac{\p\xi_1}{\p t}J
%\end{smallmatrix}
\end{pmatrix},\qquad \bm{A}_2=\begin{pmatrix}
%\begin{smallmatrix}
\frac{\p\xi_2}{\p t}J & \frac{\p\xi_2}{\p x_1}J & \frac{\p\xi_2}{\p x_2}J\\
\frac{\p\xi_2}{\p x_1}J&\frac{\p\xi_2}{\p t}J&0\\
\frac{\p\xi_2}{\p x_2}J&0&\frac{\p\xi_2}{\p t}J
%\end{smallmatrix}
\end{pmatrix}.$$
%We further denote
%$$ Q=\renewcommand\arraystretch{0.7}\begin{bmatrix}
%p\\u\\v
%\end{bmatrix},$$
Then, (\ref{eq:ale-matrix}) can be simply rewritten as
\begin{align}
\frac{d\bm{q}J}{d\tau}+\frac{\p}{\p\xi_1}\left(\bm{A}_1\bm{q}\right)+\frac{\p}{\p\xi_2}\left(\bm{A}_2\bm{q}\right)=0.
\label{eq:ale-matrix-1}
\end{align}

In Section \ref{sec:aledg}, we will derive the ALE-DG formulation for the two-dimensional acoustic wave equation. In the next section, we discuss the energy stability of a skew-symmetric ALE-DG formulation for the constant solution on moving meshes.

\section{A skew-symmetric ALE-DG formulation for mesh motion}\label{sec:mesh-ale}
In this section, we use a simple case, i.e., a constant solution on a moving mesh to illustrate energy conservation for a skew-symmetric ALE-DG formulation on moving meshes. Consider the equation
\begin{equation}
\frac{d u}{d t}=0.
\label{eq:constant}
\end{equation}
In the ALE formulation, we transfer (\ref{eq:constant}) from the time-dependent physical domain $\bm{x}=(x_1,x_2)$ with time variable $t$ onto a fixed reference domain 
${\wh{\bm{x}}}=(\xi_1,\xi_2)$ with time variable $\tau$. By setting $\bm{f}=0$ in (\ref{eq:ale4}), we obtain the corresponding system on the reference domain
\begin{equation}
\begin{split}
\frac{\p uJ}{\p \tau}+\widehat{\nabla}\cdot\left(uJ\wh{\bm{x}}_t\right)&=0,\\
\frac{\p J}{\p \tau}+\widehat{\nabla}\cdot\left(J\wh{\bm{x}}_t\right)&=0.
\end{split}
\label{eq:yamaleev}
\end{equation}
Clearly, since ${du}/{dt}=0$,
$$\qquad\frac{d}{dt}||u||^2=0,$$ and the system is energy conservative. We replicate this using a DG discretization for the reference system (\ref{eq:yamaleev}). Multiplying by test functions $v,w$ and integrating over $D^k$, we have
\begin{equation}
\begin{split}
\left(\frac{\p uJ}{\p\tau},v\right)_{L^2(D^k)}+\left(\widehat{\nabla}\cdot\left(uJ\wh{x}_t\right),v\right)_{L^2(D^k)}&=0,\\
\left(\frac{\p J}{\p\tau},w\right)_{L^2(D^k)}+\left(\widehat{\nabla}\cdot\left(J\wh{x}_t\right),w\right)_{L^2(D^k)}&=0.
\label{eq:dg-ale}
\end{split}
\end{equation}
Integrating the spatial derivative in the first equation by parts twice and using the central flux $u^*=\avg{u}$, we obtain the strong DG formulation
\begin{equation}
\begin{split}
\left(\frac{\p \left(uJ\right)}{\p\tau},v\right)_{L^2(D^k)}+\left(\widehat{\nabla}\cdot\left(uJ\wh{\bm{x}}_t\right),v\right)_{L^2(D^k)}+\frac{1}{2}\left\langle\jump{u} J\wh{\bm{x}}_t\cdot \bm{n},v\right\rangle_{L^2(\p D^k)}&=0,\\
\left(\frac{\p J}{\p\tau},w\right)_{L^2(D^k)}+\left(\widehat{\nabla}\cdot\left(J\wh{\bm{x}}_t\right),w\right)_{L^2(D^k)}&=0.
\label{eq:dg-ale-strong}
\end{split}
\end{equation}
where $\bm{n}$ is the outward normal vector on the reference element.
To show energy conservation, we use a simple trick to rewrite the first equation in (\ref{eq:dg-ale-strong}) into a skew-symmetric form
\begin{align*}
&\left(\frac{\p \left(uJ\right)}{\p \tau},v\right)_{L^2(D^k)} +\frac{1}{2} \left\langle \jump{u}J\wh{\bm{x}}_t\cdot \bm{n},v\right\rangle_{L^2(\p D^k)}\\
+&\frac{1}{2}\Big\{\left(\widehat{\nabla}\cdot\left(uJ\wh{\bm{x}}_t\right),v\right)_{L^2(D^k)}+\left(\wh{\nabla}\cdot\left(J\wh{\bm{x}}_t\right)u,v\right)_{L^2(D^k)}+\left(J\wh{\bm{x}}_t\cdot\wh{\nabla} u,v\right)_{L^2(D^k)}\Big\}
=0.
\end{align*}
Integrating the second volume term by parts yields
\begin{align}
\left(\wh{\nabla}\cdot \left(uJ\wh{\bm{x}}_t\right),v\right)_{L^2(D^k)} = \langle  uJ\wh{\bm{x}}_t\cdot \bm{n}, v\rangle_{L^2(\p D^k)}-\left(uJ\wh{\bm{x}}_t,\wh{\nabla}v\right)_{L^2(D^k)}.
\label{eq:integrate-by-part}
\end{align}
The DG formulation on the reference domain becomes
\begin{equation}
\begin{split}
\left(\frac{\p \left(uJ\right)}{\p \tau},v\right)_{L^2(D^k)}+\frac{1}{2}\langle u^+J\wh{\bm{x}}_t\cdot \bm{n}, v\rangle_{L^2(\p D^k)}\ \ \ \ \ \ \ \ \ \ \ \ \ \ \ \ \ \ \ \ \ \ \ \ \ \ \ \ \ \ \ \ \ \ \ \ \ \ &\\
+\frac{1}{2}\Big\{\left(\wh{\nabla}\cdot\left(J\wh{\bm{x}}_t\right)u,v\right)_{L^2(D^k)}+\left(J\wh{\bm{x}}_t\cdot\wh{\nabla} u,v\right)_{L^2(D^k)}-\left(uJ\wh{\bm{x}}_t,\wh{\nabla}v\right)_{L^2(D^k)}\Big\}
&=0,\\
\left(\frac{\p J}{\p\tau},w\right)_{L^2(D^k)}+\left(\widehat{\nabla}\cdot\left(J\wh{\bm{x}}_t\right),w\right)_{L^2(D^k)}&=0.
\label{eq:skew-mesh}
\end{split}
\end{equation}
In numerical experiments, the term $\wh{\nabla}\cdot\left(J\wh{\bm{x}}_t\right)$ in (\ref{eq:skew-mesh}) is computed by applying the divergence operator to the $L^2$ projection of non-polynomial function $J\wh{\bm{x}}_t$, since $J\wh{\bm{x}}_t$ may not be a polynomial of degree $N$. Note that the evolved variables are $uJ$ and $J$. It remains to specify how $u$ is determined given $uJ$ and $J$, which is the primary difference between standard DG and weight-adjusted DG.  In the following sections, we discuss energy conservation using standard DG and WADG approaches, respectively.

\subsection{Energy conservation using DG methods}\label{sec:energy-DG}
We start with energy conservation for (\ref{eq:skew-mesh}) using the DG method where weighted mass matrices are used. In the standard DG method, $u$ is defined as the solution to 
$$\left(u,vJ\right)_{L^2(D^k)}=\left(\left(uJ\right),v\right)_{L^2(D^k)},\qquad \forall v\in V_h\left(D^k\right)$$
which is equivalent to 
$$\bm{M}^k_J\bm{u}=\bm{M}^k\left(\bm{uJ}\right),\qquad \bm{u} = \left(\bm{M}^k_J\right)^{-1}\bm{M}^k\left(\bm{uJ}\right).$$
Here $\bm{u},\bm{uJ}$ are the expansion coefficients of polynomial $u$ and $uJ$. The weighted mass matrix $\bm{M}^k_J$ and the mass matrix $\bm{M}^k$ on element $\wh{D}$ are given by
$$\bm{M}^k_J=\int_{D^k}\phi^k_j\phi^k_iJ=J^k\int_{\wh{D}}\phi_j\phi_iJ=J^k\bm{M}_J,\qquad \bm{M}^k=\int_{D^k}\phi^k_j\phi^k_i=J^k\int_{\wh{D}}\phi_j\phi_iJ=J^k\bm{M},$$ 
where $\{\phi^k_i\},\{\phi_i\}$ are basis functions spanning $V_h\left(D^k\right)$ and $V_h\left(\wh{D}\right)$, respectively. The constant $J^k$ is the Jacobian of affine mapping $\bm{\Phi}^k$ between the reference element and element $D^k$. In practice,  we will use the matrix $\bm{M}_J$ and $\bm{M}$, which are referred to weighted mass matrix and mass matrix on $\wh{D}$, such that
$$\bm{u} = \left(\bm{M}^k_J\right)^{-1}\bm{M}^k\left(\bm{uJ}\right)=\left(\bm{M}_J\right)^{-1}\bm{M}\left(\bm{uJ}\right).$$

In order to show energy conservation, we take 
$$v=u,\qquad w=\frac{1}{2}\Pi_N\left(u^2\right).$$
Note that, since $\frac{\p J}{\p t},\wh{\nabla}\cdot\left(J\wh{\bm{x}}_t\right)\in P^N$, we have
\begin{align*}
\left(\frac{\p J}{\p t},\frac{1}{2}\Pi_N\left(u^2\right)\right)_{L^2(D^k)}&=\left(\frac{\p J}{\p t},\frac{1}{2}u^2\right)_{L^2(D^k)},\\
\left(\wh{\nabla}\cdot\left(J\wh{\bm{x}}_t\right),\frac{1}{2}\Pi_N\left(u^2\right)\right)_{L^2(D^k)}&=\left(\wh{\nabla}\cdot\left(J\wh{\bm{x}}_t\right),\frac{1}{2}u^2\right)_{L^2(D^k)}.
\end{align*} 
Then, summing (\ref{eq:skew-mesh}) over all elements
\begin{equation}
\begin{split}
\sum_{k}\left(\frac{\p \left(uJ\right)}{\p \tau},u\right)_{L^2(D^k)}+\sum_{k}\frac{1}{2}\left(\wh{\nabla}\cdot\left(J\wh{\bm{x}}_t\right)u,u\right)_{L^2(D^k)}\ \ \ \ \ \ \ \ \ \ \ \ \ \ \ \ \ \ \ \ \ \ \ \ \ \ \ \ \ \ \ &
\\+\sum_{k}\frac{1}{2}\Big\{\left(J\wh{\bm{x}}_t\cdot\widehat{\nabla}u,u\right)_{L^2(D^k)}-\left(uJ\wh{\bm{x}}_t,\wh{\nabla}u\right)_{L^2(D^k)}+\langle u^+J\wh{\bm{x}}_t\cdot \bm{n}, v \rangle_{L^2(\ D^k)}\Big\}&=0,\\
\sum_{k}\left(\frac{\p J}{\p\tau},\frac{u^2}{2}\right)_{L^2(D^k)}+\sum_{k}\left(\widehat{\nabla}\cdot\left(J\wh{\bm{x}}_t\right),\frac{u^2}{2}\right)_{L^2(D^k)}&=0,
\label{eq:skew-mesh1}
\end{split}
\end{equation}
For simplicity, we define the skew-symmetric term $S\left(u,v\right)$ 
$$S\left(u,v\right)=\sum_{k}\frac{1}{2}\Big\{\left(J\wh{\bm{x}}_t\cdot\widehat{\nabla}u,v\right)_{L^2(D^k)}-\left(uJ\wh{\bm{x}}_t,\wh{\nabla}v\right)_{L^2(D^k)}+\langle u^+J\wh{\bm{x}}_t\cdot \bm{n}, v \rangle_{L^2(\ D^k)}\Big\},$$
such that the first equation in (\ref{eq:skew-mesh1}) becomes
\begin{align*}
\sum_{k}\left(\frac{\p \left(uJ\right)}{\p \tau},u\right)_{L^2(D^k)}+\sum_{k}\frac{1}{2}\left(\wh{\nabla}\cdot\left(J\wh{\bm{x}}_t\right)u,u\right)_{L^2(D^k)}+S(u,u)=0.
\label{eq:skew-mesh-sum}
\end{align*}
When taking $v=u$, we immediately find that
$$\left(J\wh{\bm{x}}_t\cdot\widehat{\nabla}u,u\right)_{L^2(D^k)}-\left(uJ\wh{\bm{x}}_t,\wh{\nabla}u\right)_{L^2(D^k)}=0.$$ 
The surface contribution from the boundary is zero given reflective boundary conditions, i.e., $u=0$. Moreover, for any element-element interface, the sum of surface contributions from both neighboring physical elements is
$$\langle u^+J\wh{\bm{x}}_t\cdot \bm{n}, u \rangle_{L^2(\p D^k)}+\langle uJ\wh{\bm{x}}_t\cdot \left(-\bm{n}\right), u^+ \rangle_{L^2(\p D^k)}=0.$$
Therefore, when choosing $v=u$, the skew-symmetric term vanishes $S\left(u,u\right)=0$.

Note that, assuming continuity in time, 
%$$\int\frac{\p \left(uJ\right)}{\p\tau}u-\frac{1}{2}\frac{\p J}{\p\tau}u^2=\int\frac{\p u}{\p\tau}uJ+\frac{\p J}{\p\tau}u^2-\frac{1}{2}\frac{\p J}{\p\tau}u^2=\int\frac{\p u}{\p\tau}uJ+\frac{1}{2}\frac{\p J}{\p\tau}u^2.$$
\begin{align*}
&\left(\frac{\p\left(uJ\right)}{\p\tau},u\right)_{L^2(D^k)}-\frac{1}{2}\left(\frac{\p J}{\p\tau},u^2\right)_{L^2(D^k)}
\\=&\frac{\p}{\p\tau}\left(\left(uJ\right),u\right)_{L^2(D^k)}-\left(\left(uJ\right),\frac{\p u}{\p\tau}\right)_{L^2(D^k)}-\frac{1}{2}\left(\frac{\p J}{\p\tau},u^2\right)_{L^2(D^k)}\\
=&\frac{\p}{\p\tau}\left(u,uJ\right)_{L^2(D^k)}-\left(u,\frac{\p u}{\p\tau}J\right)_{L^2(D^k)}-\frac{1}{2}\left(\frac{\p J}{\p\tau},u^2\right)_{L^2(D^k)}\\
=&\left(\frac{\p u}{\p\tau},uJ\right)_{L^2(D^k)}+\frac{1}{2}\left(\frac{\p J}{\p\tau},u^2\right)_{L^2(D^k)}.
\end{align*}
Therefore,
$$\frac{1}{2}\frac{\p}{\p\tau}||u||_J^2=\frac{1}{2}\frac{\p}{\p\tau}\left(u,uJ\right)_{L^2\left(\wh{\Omega}_h\right)}=\left(\frac{\p u}{\p\tau},uJ\right)_{L^2\left(\wh{\Omega}_h\right)}+\frac{1}{2}\left(\frac{\p J}{\p\tau},u^2\right)_{L^2\left(\wh{\Omega}_h\right)}.$$
Subtracting the second equation from the first equation in (\ref{eq:skew-mesh1}), we obtain
\begin{equation}
\frac{1}{2}\frac{\p}{\p\tau}||u||^2_J+S\left(u,u\right)+\sum_{k}\frac{1}{2}\left(\left(\wh{\nabla}\cdot\left(J\wh{\bm{x}}_t\right)u,u\right)_{L^2(D^k)}-\left(\widehat{\nabla}\cdot\left(J\wh{\bm{x}}_t\right),u^2\right)_{L^2(D^k)}\right)=0.
\label{eq:skew-energy}
\end{equation}
By skew-symmetry of $S\left(u,v\right)$ and the fact that $$\left(\wh{\nabla}\cdot\left(J\wh{\bm{x}}_t\right)u,u\right)_{L^2(D^k)}=\left(\wh{\nabla}\cdot\left(J\wh{\bm{x}}_t\right),u^2\right)_{L^2(D^k)},$$
when integrals are computed using quadrature, we have the following result.
\begin{thm}\label{thm:energy-conservation}
	The skew-symmetric formulation (\ref{eq:skew-mesh}) using the standard DG method is energy conservative in the sense that
	\begin{equation}
	\frac{1}{2}\frac{\p}{\p\tau}||u||^2_J=0.
	\label{eq:skew-u}
	\end{equation}
\end{thm}

\subsection{Weight-adjusted DG methods}
Note that, when using the standard DG approach, we encounter the matrix $\bm{M}_J$ when computing original solution $u$ from evolving variables $uJ$ and $J$. At each time step, the inversion of the reference weighted mass matrices $\bm{M}_J$ is required. When $J$ is approximated by a constant over each element, we are able to apply $\left(\bm{M}_J\right)^{-1}$ for all elements using a single reference mass matrix inverse $\bm{M}^{-1}$ with scale factor $1/J$. By definition, we have
$$ \left(\bm{M}_J\right)^{-1}=\frac{1}{J}\bm{M}^{-1}.$$
However, when $J$ possesses sub-element variations, inverses of weighted mass matrices become distinct from element to element. Typical implementations require precomputation and storage of these weighted mass matrix inverses  \cite{mercerat2015nodal, bencomo2015discontinuous}, which significantly increases the storage cost of high order schemes.

To address this issue, we approximate the weighted mass matrix using an easily-invertible weight-adjusted approximation \cite{chan2017weight,chan2018weight,chan2017curved}, which is energy stable and high order accurate for sufficiently regular weighting functions. The weight-adjusted approximation $\tilde{\bm{M}}_J$ is given by
$$\bm{M}_J\approx \tilde{\bm{M}}_J =\bm{M}\left(\bm{M}_{1/J}\right)^{-1}\bm{M}.$$
The inverse of the approximation is then
$$\left(\bm{M}_J\right)^{-1}\approx \left(\tilde{\bm{M}}_J\right)^{-1}=\bm{M}^{-1}\bm{M}_{1/J}\bm{M}^{-1}.$$
Note that the weighting function now only appear in $M_{1/J}$, which can be applied in a low-storage quadrature-based manner using sufficiently accurate quadrature rules. We choose simplicial quadratures that are exact for polynomials of degree $2N+1$ \cite{xiao2010numerical}, and let ${\bm{r}}^q_i,{\bm{w}}_i$ denote the quadrature points and weights on the reference element $\wh{D}$. We define the reference interpolation matrix $\bm{V}_q$ as 
$$\left(\bm{V}_q\right)_{ij} = \phi_j\left({\bm{r}}^q_i\right),$$
whose columns consist of values of basis functions at quadrature points. Then,
\begin{align*}
\bm{M}=\bm{V}_q^T\textmd{diag}\left({\bm{w}}\right)\bm{V}_q,\qquad
\bm{M}_J=\bm{V}_q^T\textmd{diag}\left({\bm{w}}\right)\textmd{diag}\left(\bm{d}\right)\bm{V}_q,\qquad \bm{d}_i=J\left(X\left(\bm{\Phi}^k{\bm{r}}^q_i,\tau\right)\right)
\end{align*}
where $X\left(\bm{\Phi}^k{\bm{x}}^q_i,\tau\right)$ are quadrature points on the physical element at time $\tau$, and $\bm{d}$ denote the values of the Jacobian at these quadrature points. Using the weight-adjusted approximation, we have
$$\bm{u}=\left(\bm{M}_J\right)^{-1}\bm{M}\left(\bm{uJ}\right)\approx\bm{M}^{-1}\bm{M}_{1/J}\left(\bm{uJ}\right).$$
We need to apply the product of an unweighted mass matrix inverse and a weighted mass matrix, which can be done using quadrature-based matrices as follows: 
\begin{align}
\bm{M}^{-1}\bm{M}_{1/J}=\bm{P}_q\textmd{diag}\left(\bm{d}\right)\bm{V}_q,
\label{eq:wadg}
\end{align}
where $\bm{P}_q=\bm{M}^{-1}\bm{V}_q^T\textmd{diag}\left({\bm{w}}\right)$ is a quadrature discretization of the polynomial $L^2$ projection operator on the reference element. Since $\bm{P}_q$ and $\bm{V}_q$ are reference operators, the application of the weight-adjusted approximation only require $O(N^2)$ storage for values of the Jacobian at quadrature points for each element. In contrast, using standard DG approach, storing the weighted mass matrix inverses or factorizations requires $O(N^4)$ storage on each element. For a general $d$ dimensional element, the cost of matrix assembly using quadrature and solving the resulting matrix system are both $O(N^{3d})$ over each element.
\subsection{Energy conservation using WADG methods}\label{sec:energy-wadg}
The proof of energy conservation for WADG methods is slightly different from the proof of energy conservation for DG methods. Using the weight-adjusted approximation to the weighted mass matrix, we can derive an upper bound for the energy variation to show an asymptotic discrete energy conservation property. Note that unlike standard DG, the energy for WADG is measured in the norm $||\cdot||_{1/J}$, which is defined as
$$E=||uJ||^2_{1/J}=\left(\frac{\left(uJ\right)}{J},\left(uJ\right)\right)_{L^2\left(\wh{\Omega}\right)}.$$ 

In numerical implementations, the field variables are $uJ$ and the Jacobian $J$, which are evolved at each time step. In WADG, we replace weighted mass matrices by weight-adjusted mass matrices 
$$\bm{M}_J\approx \bm{M}\left(\bm{M}_{1/J}\right)^{-1}\bm{M},$$
and solve $\bm{u}$ through
$$\bm{u}=\bm{M}^{-1}\bm{M}_{1/J}\left(\bm{uJ}\right),$$
which is equivalent to defining $u$ as the solution to 
$$\left(u,v\right)_{L^2(D^k)}=\left(\frac{\left(uJ\right)}{J},v\right)_{L^2(D^k)},\qquad\forall v\in V_h\left(D^k\right).$$
Thus, using WADG, we have
$$u=\Pi_N\tilde{u},\qquad \tilde{u}=\frac{\left(uJ\right)}{J},$$
where $uJ,J\in P^N$ and $\Pi_N$ is the $L^2$ projection operator onto degree $N$ polynomials. In this section, we only use properties of the $L^2$ projection operator which hold when integrals are approximated using quadrature.

To show discrete energy conservation, we follow the way used in the previous section with slight modifications. In (\ref{eq:skew-mesh1}), we choose
$$v=u,\qquad w=\frac{1}{2}\Pi_N\left(\tilde{u}^2\right),$$
and sum over elements
\begin{equation}
\begin{split}
\sum_{k}\left(\frac{\p \left(uJ\right)}{\p\tau},u\right)_{L^2(D^k)}+S\left(u,u\right)+\sum_{k}\frac{1}{2}\left(\wh{\nabla}\cdot\left(J\wh{\bm{x}}_t\right)u,u\right)_{L^2(D^k)}&=0,\\
\sum_{k}\frac{1}{2}\left(\frac{\p J}{\p \tau},\Pi_N\left(\tilde{u}^2\right)\right)_{L^2(D^k)}+\sum_{k}\frac{1}{2}\left(\wh{\nabla}\cdot\left(J\wh{\bm{x}}_t\right),\Pi_N\left(\tilde{u}^2\right)\right)_{L^2(D^k)}&=0.
\label{eq:dis-skew1}
\end{split}
\end{equation}
Again, by skew-symmetry of $S(u,v)$,
$$S\left(u,u\right)=0.$$
Now, we subtract the second equation in (\ref{eq:dis-skew1}) from the first equation
\begin{equation}
\begin{split}
&\sum_{k}\left(\frac{\p \left(uJ\right)}{\p\tau},u\right)_{L^2(D^k)}-\sum_{k}\frac{1}{2}\left(\frac{\p J}{\p \tau},\Pi_N\left(\tilde{u}^2\right)\right)_{L^2(D^k)}\\+&\sum_{k}\frac{1}{2}\left(\wh{\nabla}\cdot\left(J\wh{\bm{x}}_t\right)u,u\right)_{L^2(D^k)}
-\sum_{k}\frac{1}{2}\left(\wh{\nabla}\cdot\left(J\wh{\bm{x}}_t\right),\Pi_N\left(\tilde{u}^2\right)\right)_{L^2(D^k)}=0.
\end{split}
\label{eq:dis-skew2}
\end{equation}
Note that
%$$\textcolor{black}{\frac{\p}{\p \tau}||uJ||_{1/J}^2}=\frac{\p}{\p \tau}\int_{\Omega_h}\frac{\left(uJ\right)^2}{J}=\int_{\Omega_h}\frac{\p}{\p \tau}\frac{\left(uJ\right)^2}{J}=2\int_{\Omega_h}\frac{\left(uJ\right)}{J}\frac{\p uJ}{\p \tau}-\int_{\Omega_h}\frac{\p J}{\p \tau}\frac{\left(uJ\right)^2}{J^2}.$$
\begin{align*}
\frac{\p}{\p \tau}||uJ||_{1/J}^2&=\sum_{k}\frac{\p}{\p \tau}\left(\frac{\left(uJ\right)}{J},\left(uJ\right)\right)_{L^2\left(D^k\right)}\\
&=\sum_{k}\left(\frac{\p}{\p\tau}\left(\frac{\left(uJ\right)}{J}\right), \left(uJ\right)\right)_{L^2\left(D^k\right)}+\sum_{k}\left(\frac{\left(uJ\right)}{J}, \frac{\p\left(uJ\right)}{\p\tau}\right)_{L^2\left(D^k\right)}\\
&=2\sum_{k}\left(\frac{\left(uJ\right)}{J}, \frac{\p\left(uJ\right)}{\p\tau}\right)_{L^2\left(D^k\right)}-\sum_{k}\left(\frac{\p J}{\p\tau}, \frac{\left(uJ\right)^2}{J^2}\right)_{L^2\left(D^k\right)}\\
&=2\sum_{k}\left(\tilde{u}, \frac{\p\left(uJ\right)}{\p\tau}\right)_{L^2\left(D^k\right)}-\sum_{k}\left(\frac{\p J}{\p\tau}, \left(\tilde{u}\right)^2\right)_{L^2\left(D^k\right)}.
\end{align*}
Since $uJ$ and $J$ are degree $N$ polynomials in space on $D^k$, $\frac{\p uJ}{\p t},\frac{\p J}{\p t}\in P^N$, and we can use properties of the $L^2$ projection operator to show
\begin{align*}
\frac{\p}{\p \tau}||uJ||_{1/J}^2&=2\sum_{k}\left(\tilde{u}, \frac{\p\left(uJ\right)}{\p\tau}\right)_{L^2\left(D^k\right)}-\sum_{k}\left(\frac{\p J}{\p\tau}, \left(\tilde{u}\right)^2\right)_{L^2\left(D^k\right)}\\&=2\sum_{k}\left(u, \frac{\p\left(uJ\right)}{\p\tau}\right)_{L^2\left(D^k\right)}-\sum_{k}\left(\frac{\p J}{\p\tau}, \Pi_N\left(\tilde{u}\right)^2\right)_{L^2\left(D^k\right)}.
\end{align*}
Therefore, from (\ref{eq:dis-skew2}), we have
\begin{equation*}
\begin{split}
\frac{\p}{\p \tau}||uJ||_{1/J}^2+\sum_{k}\left(\wh{\nabla}\cdot\left(J\wh{\bm{x}}_t\right)u,u\right)_{L^2\left(D^k\right)}
-\sum_{k}\left(\wh{\nabla}\cdot\left(J\wh{\bm{x}}_t\right),\Pi_N\left(\tilde{u}^2\right)\right)_{L^2\left(D^k\right)}=0.
\end{split}
\end{equation*}
For the last term, we have
$$\left(\wh{\nabla}\cdot\left(J\wh{\bm{x}}_t\right),\Pi_N\left(\tilde{u}^2\right)\right)_{L^2\left(D^k\right)}=\left(\wh{\nabla}\cdot\left(J\wh{\bm{x}}_t\right),\tilde{u}^2\right)_{L^2\left(D^k\right)}=\left(\wh{\nabla}\cdot\left(J\wh{\bm{x}}_t\right)\tilde{u},\tilde{u}\right)_{L^2\left(D^k\right)},$$
and
\begin{align*}
\frac{\p}{\p \tau}||uJ||_{1/J}^2+\sum_{k}\left(\wh{\nabla}\cdot\left(J\wh{\bm{x}}_t\right)u,u\right)_{L^2\left(D^k\right)}
-\sum_{k}\left(\wh{\nabla}\cdot\left(J\wh{\bm{x}}_t\right)\tilde{u},\tilde{u}\right)_{L^2\left(D^k\right)}=0.
\end{align*}
Integrating over time, we obtain
\begin{equation*}
\begin{split}
\int_0^T\textcolor{black}{\frac{\p}{\p \tau}||uJ||_{1/J}^2}+\sum_{k}\int_0^T\left(\wh{\nabla}\cdot\left(J\wh{\bm{x}}_t\right)u,u\right)_{L^2\left(D^k\right)}
-\sum_{k}\int_0^T\left(\wh{\nabla}\cdot\left(J\wh{\bm{x}}_t\right)\tilde{u},\tilde{u}\right)_{L^2\left(D^k\right)}=0.
\end{split}
\end{equation*}
Therefore, 
\begin{equation}
\begin{split}
\Big|||uJ\left(\cdot,T\right)||_{1/J}^2-||uJ\left(\cdot,0\right)||_{1/J}^2\Big|&\leq\sum_{k} \int_0^T\Big|\left(\wh{\nabla}\cdot\left(J\wh{\bm{x}}_t\right)u,u\right)_{L^2\left(D^k\right)}-\left(\wh{\nabla}\cdot\left(J\wh{\bm{x}}_t\right)\tilde{u},\tilde{u}\right)_{L^2\left(D^k\right)}\Big|\\
&\leq\max_{t\in[0,T]}||\wh{\nabla}\cdot\left(J\wh{\bm{x}}_t\right)||_{\infty}\sum_k\int_{0}^{T} \Big|\left(u,u\right)_{L^2\left(D^k\right)}-\left(\tilde{u},\tilde{u}\right)_{L^2\left(D^k\right)}\Big| \\
&\leq \max_{t\in[0,T]}||\wh{\nabla}\cdot\left(J\wh{\bm{x}}_t\right)||_{\infty}\sum_k\int_{0}^{T}\Big|\left(\tilde{u},u\right)_{L^2\left(D^k\right)}-\left(\tilde{u},\tilde{u}\right)_{L^2\left(D^k\right)}\Big|\\
&\leq \max_{t\in[0,T]}||\wh{\nabla}\cdot\left(J\wh{\bm{x}}_t\right)||_{\infty}\sum_k\int_{0}^{T}\Big|\left(\tilde{u},u-\tilde{u}\right)_{L^2\left(D^k\right)}\Big|.
\end{split}
\label{eq:dis-skew3}
\end{equation}
Since $u=\Pi_N\tilde{u}$, we have
$$\left(u,u-\tilde{u}\right)_{L^2\left(D^k\right)}=\left(u,\Pi_N\tilde{u}-\tilde{u}\right)_{L^2\left(D^k\right)}=0,$$
by the fact that $\Pi_N$ is an orthogonal projection operator.
Thus, from (\ref{eq:dis-skew3}), we can derive that
\begin{align*}
\Big|||uJ\left(\cdot,T\right)||_{1/J}^2-||uJ\left(\cdot,0\right)||_{1/J}^2\Big|&\leq \max_{t\in[0,T]}||\wh{\nabla}\cdot\left(J\wh{\bm{x}}_t\right)||_{\infty}\sum_k\int_{0}^{T}\Big|\left(\tilde{u},u-\tilde{u}\right)_{L^2\left(D^k\right)}\Big|\\
&\leq \max_{t\in[0,T]}||\wh{\nabla}\cdot\left(J\wh{\bm{x}}_t\right)||_{\infty}\sum_k\int_{0}^{T}\Big|\left(\tilde{u}-u,u-\tilde{u}\right)_{L^2\left(D^k\right)}\Big|\\
&\leq \max_{t\in[0,T]}||\wh{\nabla}\cdot\left(J\wh{\bm{x}}_t\right)||_{\infty}\sum_k\int_{0}^{T}||\tilde{u}-\Pi_N\tilde{u}||_{L^2\left(D^k\right)}^2\\
&\leq \max_{t\in[0,T]}||\wh{\nabla}\cdot\left(J\wh{\bm{x}}_t\right)||_{\infty}\sum_k\int_0^TC(t)h^{2N+2}||u||^2_{H^{N+1}\left(D^k\right)}\\
&\leq \max_{t\in[0,T]}||\wh{\nabla}\cdot\left(J\wh{\bm{x}}_t\right)||_{\infty}\int_0^TC(t)h^{2N+2}||u||^2_{H^{N+1}\left(\wh{\Omega}_h\right)}.
\end{align*}
For the last inequality, we use the standard interpolation estimate for sufficiently regular functions \cite{chan2017weight}
$$||\tilde{u}-\Pi_N\tilde{u}||^2\leq C(t)h^{2N+2}||u||^2_{H^{N+1}\left(D^k\right)},$$
where $||\cdot||_{H^{N+1}}$ denotes the $L^2$ Sobolev norm of degree $N+1$. Assuming $C(t)$ is bounded such that $C(t)\leq C_{\max}$ for $t\in[0,T]$ (which is true if the regularity of $\widetilde{u}$ does not change for $t\in[0,T]$), then
\begin{align*}
\Big|||uJ\left(\cdot,T\right)||_{1/J}^2-||uJ\left(\cdot,0\right)||_{1/J}^2\Big|&\leq \max_{t\in[0,T]}||\wh{\nabla}\cdot\left(J\wh{\bm{x}}_t\right)||_{\infty}\max_{t\in[0,T]}||u||^2_{H^{N+1}\left(\wh{\Omega}_h\right)}C_{\max}\int^T_0h^{2N+2}\\
&\leq \max_{t\in[0,T]}||\wh{\nabla}\cdot\left(J\wh{\bm{x}}_t\right)||_{\infty}\max_{t\in[0,T]}||u||^2_{H^{N+1}\left(\wh{\Omega}_h\right)}C_{\max}Th^{2N+2}.
\end{align*}
Therefore, we obtain the following theorem.
\begin{thm}\label{thm:upper-bound}
	The skew-symmetric formulation (\ref{eq:skew-mesh}) using the WADG method has an upper bound for the energy variation given by
	\begin{equation}
	\textcolor{black}{\Big|||uJ\left(\cdot,T\right)||_{1/J}^2-||uJ\left(\cdot,0\right)||_{1/J}^2\Big|\leq Ch^{2N+2},}
	\label{eq:estimate}
	\end{equation}
	for fixed $T$ and sufficiently regular solution $u(x,t)$.
\end{thm}
Theorem \ref{thm:upper-bound} implies that, for sufficiently regular solutions, the skew-symmetric DG discretization (\ref{eq:skew-mesh}) is asymptotically energy conservative as we shrink the mesh size or increase the order of approximation. Although this is not strictly energy conservative, the result is sufficient to ensure energy boundedness and stability in practice. In Section \ref{sec:ale-experiment}, we verify the validation of this upper bound using numerical experiments. In the next section, I will use the similar techniques to derive a skew-symmetric DG formulation for the acoustic wave equation on moving meshes.

\subsection{Energy equivalence using both standard DG and WADG approaches}
In Section~\ref{sec:energy-DG} and \ref{sec:energy-wadg}, we discussed energy conservation for the ALE-DG formulation (\ref{eq:skew-mesh}) using standard DG and weight-adjusted DG approach. Note that we use two different energy norms $||\cdot||_J$ and $||\cdot||_{1/J}$ to measure energy variation. In this part, we show that the two measures of energy $||u||_J$ and $||uJ||_{1/J}$ are equivalent when using weight-adjusted DG approach, such that
$$c_1||u||_J\leq ||uJ||_{1/J}\leq c_2||u||_J,$$
where $u$ is defined through
$$\left(u,v\right)_{L^2(D^k)}=\left(\frac{\left(uJ\right)}{J},v\right)_{L^2(D^k)},\qquad\forall v\in V_h\left(D^k\right).$$ 
We first assume $J$ is bounded over $\wh{\Omega}_h$ by
$$0<J_{\min}\leq J\leq J_{\max}.$$
By definition,
$$||uJ||^2_{1/J}=\sum_k\left(\frac{\left(uJ\right)}{J},\left(uJ\right)\right)_{L^2(D^k)}=\sum_k\left(\Pi_N\tilde{u},\left(uJ\right)\right)_{L^2(D^k)}=\sum_k\left(u,\left(uJ\right)\right)_{L^2(D^k)}.$$
On element $D^k$, we multiply and divide by $J$ to get
\begin{align*}
\left(u,\left(uJ\right)\right)_{L^2(D^k)}=\left(u,\frac{\left(uJ\right)}{J}J\right)_{L^2(D^k)}\leq  J_{\max}\left(u,\frac{\left(uJ\right)}{J}\right)_{L^2(D^k)}=J_{\max}\left(u,\Pi_N\frac{\left(uJ\right)}{J}\right)_{L^2(D^k)}.
\end{align*}
Thus,
$$\left(u,\left(uJ\right)\right)_{L^2(D^k)}\leq J_{\max}\left(u,u\right)_{L^2(D^k)}.$$
Multiplying and dividing by $J$ again gives
$$\left(u,\left(uJ\right)\right)_{L^2(D^k)}\leq J_{\max}\left(u,u\right)_{L^2(D^k)}=J_{\max}\left(u,\frac{uJ}{J}\right)_{L^2(D^k)}\leq \frac{J_{\max}}{J_{\min}}\left(u,uJ\right)_{L^2(D^k)}.$$
Therefore,
$$||uJ||^2_{1/J}\leq \sum_k\frac{J_{\max}}{J_{\min}}\left(u,uJ\right)_{L^2(D^k)}=\frac{J_{\max}}{J_{\min}}||u||_J^2.$$
On the other hand, 
\begin{align*}
\left(u,\left(uJ\right)\right)_{L^2(D^k)}=\left(u,\frac{\left(uJ\right)}{J}J\right)_{L^2(D^k)}\geq  J_{\min}\left(u,\frac{\left(uJ\right)}{J}\right)_{L^2(D^k)}=J_{\min}\left(u,u\right)_{L^2(D^k)},
\end{align*}
and
$$\left(u,\left(uJ\right)\right)_{L^2(D^k)}\geq J_{\min}\left(u,u\right)_{L^2(D^k)}=J_{\min}\left(u,\frac{uJ}{J}\right)_{L^2(D^k)}\geq \frac{J_{\min}}{J_{\max}}\left(u,uJ\right)_{L^2(D^k)}.$$
So
$$||uJ||^2_{1/J}\geq \sum_k\frac{J_{\min}}{J_{\max}}\left(u,uJ\right)_{L^2(D^k)}=\frac{J_{\min}}{J_{\max}}||u||_J^2.$$
Finally, we end up with the following theorem:
\begin{thm}
\label{thm:eqn}
Using the weight-adjusted DG method, energy measured by these two aforementioned norms are equivalent,
$$\sqrt{\frac{J_{\min}}{J_{\max}}}||u||_J\leq||uJ||_{1/J}\leq \sqrt{\frac{J_{\max}}{J_{\min}}}||u||_J,$$
where $u,uJ\in P^N$ are connected through
$$\left(u,v\right)_{L^2(D^k)}=\left(\frac{\left(uJ\right)}{J},v\right)_{L^2(D^k)},\qquad\forall v\in V_h\left(D^k\right).$$ 
\end{thm}
Note that, by Theorem~\ref{thm:upper-bound}, the energy computed by $||uJ||_{1/J}$ is bounded at time $T$. Therefore, we can immediately conclude that energy computed by $||u||_J$ is also bounded by Theorem~\ref{thm:eqn}.
\section{Novel DG methods for wave propagation on moving meshes}\label{sec:aledg}
In this section, we present an ALE-DG formulation for wave propagation on moving curved meshes. Similar to the previous section, we use a skew-symmetric formulation under which we can easily prove energy stability. We start with (\ref{eq:ale-matrix-1}) and the geometric conservation Law (GCL)
\begin{equation}
\begin{split}
\frac{d\left(\bm{q}J\right)}{d\tau}+\frac{\p}{\p\xi_1}\left(A^1\bm{q}\right)+\frac{\p}{\p\xi_2}\left(A^2\bm{q}\right)&=0,\\
\frac{\p J}{\p \tau}+\wh{\nabla}\cdot\left(J\wh{\bm{x}}_t\right)&=0,
\end{split}
\label{eq:aledg}
\end{equation}
where
$$A^1=\begin{pmatrix}
%\begin{smallmatrix}
\frac{\p\xi_1}{\p t}J & \frac{\p\xi_1}{\p x_1}J & \frac{\p\xi_1}{\p x_2}J\\
\frac{\p\xi_1}{\p x_1}J&\frac{\p\xi_1}{\p t}J&0\\
\frac{\p\xi_1}{\p x_2}J&0&\frac{\p\xi_1}{\p t}J
%\end{smallmatrix}
\end{pmatrix},\qquad A^2=\begin{pmatrix}
%\begin{smallmatrix}
\frac{\p\xi_2}{\p t}J & \frac{\p\xi_2}{\p x_1}J & \frac{\p\xi_2}{\p x_2}J\\
\frac{\p\xi_2}{\p x_1}J&\frac{\p\xi_2}{\p t}J&0\\
\frac{\p\xi_2}{\p x_2}J&0&\frac{\p\xi_2}{\p t}J
%\end{smallmatrix}
\end{pmatrix}.$$
To construct a variational formulation, we multiply the above equation by test functions $\bm{w},\theta$ and integrate over $D^k$
\begin{equation}
\begin{split}
\left(\frac{d\left(\bm{q}J\right)}{d\tau},\bm{w}\right)_{L^2(D^k)}&=-\left(\frac{\p}{\p\xi_1}\left(A^1\bm{q}\right),\bm{w}\right)_{L^2(D^k)}-\left(\frac{\p}{\p\xi_2}\left(A^2\bm{q}\right),\bm{w}\right)_{L^2(D^k)},\\
\left(\frac{\partial J}{d\tau},\theta\right)_{L^2(D^k)}&=-\left(\wh{\nabla}\cdot\left(J\wh{\bm{x}}_t\right),\theta\right)_{L^2(D^k)}.
\label{eq:acous-ale-dg1}
\end{split}
\end{equation}
Now, we rewrite (\ref{eq:acous-ale-dg1}) into a skew-symmetric form
\begin{equation}
\begin{split}
\left(\frac{d\left(\bm{q}J\right)}{d\tau},\bm{w}\right)_{L^2(D^k)}=&-\frac{1}{2}\left(A^1\left(\frac{\p}{\p\xi_1}\bm{q}\right), \bm{w}\right)_{L^2(D^k)}+\frac{1}{2}\left(A^1\bm{q},\frac{\p}{\p\xi_1}\bm{w}\right)_{L^2(D^k)}\\
&-\frac{1}{2}\left(\left(\frac{\p}{\p\xi_1}A^1\right)\bm{q},\bm{w}\right)_{L^2(D^k)}-\frac{1}{2}\left\langle A^1\bm{q}^*\widehat{n}_1,\bm{w}\right\rangle_{L^2(\p D^k)}\\
&-\frac{1}{2}\left(A^2\left(\frac{\p}{\p\xi_2}\bm{q}\right), \bm{w}\right)_{L^2(D^k)}+\frac{1}{2}\left(A^2\bm{q},\frac{\p}{\p\xi_2}\bm{w}\right)_{L^2(D^k)}\\
&-\frac{1}{2}\left(\left(\frac{\p}{\p\xi_2}A^2\right)\bm{q},\bm{w}\right)_{L^2(D^k)}-\frac{1}{2}\left\langle A^2\bm{q}^*\widehat{n}_2,\bm{w}\right\rangle_{L^2(\p D^k)}\\
\left(\frac{\partial J}{d\tau},\theta\right)_{L^2(D^k)}=&-\left(\wh{\nabla}\cdot\left(J\wh{\bm{x}}_t\right),\theta\right)_{L^2(D^k)},
\label{eq:acous-ale-dg12}
\end{split}
\end{equation}
where $\bm{q}^*$ is the numerical flux and $\widehat{n}=\left(\widehat{n}_1,\widehat{n}_2\right)$ is the reference normal vector. We combine all surface terms together and obtain
\begin{equation}
\begin{split}
\left(\frac{d\left(\bm{q}J\right)}{d\tau},\bm{w}\right)_{L^2(D^k)}=&-\frac{1}{2}\left(\frac{\p}{\p\xi_1}\left(A^1\bm{q}\right),\bm{w}\right)_{L^2(D^k)}+\frac{1}{2}\left(\bm{q},\frac{\p}{\p\xi_1}\left(A^1\bm{w}\right)\right)_{L^2(D^k)}\\&-\frac{1}{2}\left(\left(\frac{\p}{\p\xi_1}A^1\right)\bm{q},\bm{w}\right)_{L^2(D^k)}\\
&-\frac{1}{2}\left(\frac{\p}{\p\xi_2}\left(A^2\bm{q}\right),\bm{w}\right)_{L^2(D^k)}+\frac{1}{2}\left(\bm{q},\frac{\p}{\p\xi_2}\left(A^2\bm{w}\right)\right)_{L^2(D^k)}\\&-\frac{1}{2}\left(\left(\frac{\p}{\p\xi_2}A^2\right)\bm{q},\bm{w}\right)_{L^2(D^k)}\\
&-\frac{1}{2}\left\langle \bm{q}^*,A_n\bm{w}\right\rangle_{L^2(\p D^k)}\\
\left(\frac{\partial J}{d\tau},\theta\right)_{L^2(D^k)}=&-\left(\wh{\nabla}\cdot\left(J\wh{\bm{x}}_t\right),\theta\right)_{L^2(D^k)}.
\label{eq:acous-ale-dg13}
\end{split}
\end{equation}
where $A_n=A^1\widehat{n}_1+A^2\widehat{n}_2$.

The choice of numerical fluxes is important for DG methods. Traditional upwind numerical fluxes, which are based on solvers for local Riemann problems, are difficult to compute analytically, e.g., for anisotropic elastic media \cite{zhan2018exact}. There is also no guarantee of energy stability for upwind fluxes in the presence of sub-cell media heterogeneities. Instead of the upwind flux, we use the so-called penalty flux in our proposed DG scheme for wave equations on moving meshes. This penalty flux will penalize the residual of appropriate continuity conditions and add numerical dissipation to damp out spurious components of the solution. The choice of penalty fluxes is simple and straightforward, as it is constructed based only on appropriate continuity conditions. In this work, I propose a new penalty flux for the ALE-DG formulation, which is provably consistent, energy stable and has a simple form. The proposed penalty flux is motivated by the surface term in (\ref{eq:acous-ale-dg13}) and is given by
\begin{align}
\bm{q}^*=\bm{q}^+-\tau_{q} A_n\jump{\bm{q}}.
\label{eq:ale-flux}
\end{align}
where $\tau_q\geq 0$ is penalty parameter. The proposed penalty flux can be considered as an extension of the standard penalty flux \cite{Hesthaven2007} for static meshes since it will reduce to the standard penalty flux if there is no mesh motion. In the following sections, I will prove the consistency and energy stability for this choice of penalty flux. In Section \ref{sec:ale-experiment}, we present numerical experiments to verify the high order accuracy of the proposed DG method.

\subsection{Consistency}
In this part, we prove that the DG formulation (\ref{eq:acous-ale-dg13}) with numerical fluxes (\ref{eq:ale-flux}) is consistent under appropriate continuity conditions. Assume that $\bm{q}$ is exact solution of the acoustic wave equation, and that stable and consistent boundary conditions as described in \cite{chan2017weight,chan2018weight} are imposed through consistent modifications of the numerical flux. Now, we integrate the volume terms in (\ref{eq:acous-ale-dg13}) involving derivatives of test functions by parts again. Then, plugging these exact solutions into (\ref{eq:acous-ale-dg13}) causes the volume terms to vanish. Consistency follows if the numerical flux terms also vanish.  After integration by parts, the flux term become
\begin{equation*}
\begin{split}
\langle\bm{q}^*-\bm{q},A_n\bm{w}\rangle = \langle A_n\left(\bm{q}^*-\bm{q}\right),\bm{w}\rangle=\langle A_n\left(\jump{\bm{q}}-\tau_qA_n\jump{\bm{q}}\right),\bm{w}\rangle
\end{split}
\label{eq:ale-flux1}
\end{equation*}
The exact solutions across inter-element interfaces should satisfy $\jump{\bm{q}}=0$, and we can immediately conclude
$$\langle\bm{q}^*-\bm{q},A_n\bm{w}\rangle=0.$$
Therefore, all numerical flux terms vanish and the skew-symmetric ALE-DG formulation with penalty flux (\ref{eq:ale-flux}) is consistent.

\subsection{Energy stability}
The formulations (\ref{eq:acous-ale-dg13}) with numerical fluxes (\ref{eq:ale-flux}) can also be shown to be energy stable for $\tau_q\geq 0$. For simplicity, we assume zero homogeneous Dirichlet boundary conditions on $\partial \Omega$ in the proof of energy stability. Similarly, using DG methods, we take 
$$\bm{w} = \bm{q},\quad h = \frac{1}{2}\Pi_N\left(\bm{q}^T\bm{q}\right)$$ 
in formulation (\ref{eq:acous-ale-dg13}) and sum over elements. The first equation becomes
\begin{equation}
\begin{split}
\left(\frac{d\left(\bm{q}J\right)}{d\tau},\bm{q}\right)_{L^2\left(\wh{\Omega}_h\right)}=&-\frac{1}{2}\left(\left(\frac{\p}{\p\xi_1}A^1\right)\bm{q},\bm{q}\right)_{L^2\left(\wh{\Omega}_h\right)}
-\frac{1}{2}\left(\left(\frac{\p}{\p\xi_2}A^2\right)\bm{q},\bm{q}\right)_{L^2\left(\wh{\Omega}_h\right)}\\
&-\sum_{k}\frac{1}{2}\left\langle \bm{q}^*,A_n\bm{q}\right\rangle_{L^2(\p D^k)}.
\label{eq:ale-sum}
\end{split}
\end{equation}
The second equation becomes
\begin{equation}
\begin{split}
\frac{1}{2}\left(\frac{\partial J}{d\tau},\Pi_N\left(\bm{q}^T\bm{q}\right)\right)_{L^2\left(\wh{\Omega}_h\right)}&=-\frac{1}{2}\left(\wh{\nabla}\cdot\left(J\wh{\bm{x}}_t\right),\Pi_N\left(\bm{q}^T\bm{q}\right)\right)_{L^2\left(\wh{\Omega}_h\right)}.
\end{split}
\label{eq:ale-J}
\end{equation}
Again, since $\frac{\p J}{\p\tau},\wh{\nabla}\cdot\left(J\wh{\bm{x}}_t\right)\in P^N$, then (\ref{eq:ale-J}) is equivalent to
\begin{equation}
\begin{split}
\frac{1}{2}\left(\frac{\partial J}{d\tau},\bm{q}^T\bm{q}\right)_{L^2\left(\wh{\Omega}_h\right)}&=-\frac{1}{2}\left(\wh{\nabla}\cdot\left(J\wh{\bm{x}}_t\right),\bm{q}^T\bm{q}\right)_{L^2\left(\wh{\Omega}_h\right)}.
\end{split}
\label{eq:ale-J1}
\end{equation}	  
Subtracting (\ref{eq:ale-J1}) from (\ref{eq:ale-sum}) gives, on the left hand side
\begin{align*}
\left(\frac{d\left(\bm{q}J\right)}{d\tau},\bm{q}\right)_{L^2\left(\wh{\Omega}_h\right)}-\frac{1}{2}\left(\frac{\partial J}{d\tau},\bm{q}^T\bm{q}\right)_{L^2\left(\wh{\Omega}_h\right)}&=\frac{1}{2}\frac{d}{d\tau}\left(\bm{q}^T\bm{q},J\right)_{L^2\left(\wh{\Omega}_h\right)}\\&=\frac{1}{2}\frac{\p}{\p\tau}\left(||p||_J^2+||u||_J^2+||v||_J^2\right),
\end{align*}
and on the right hand side
\begin{align}
RHS=&-\frac{1}{2}\left(\left(\frac{\p}{\p\xi_1}A^1\right)\bm{q},\bm{q}\right)_{L^2\left(\wh{\Omega}_h\right)}
-\frac{1}{2}\left(\left(\frac{\p}{\p\xi_2}A^2\right)\bm{q},\bm{q}\right)_{L^2\left(\wh{\Omega}_h\right)}\notag\\&+\frac{1}{2}\left(\wh{\nabla}\cdot\left(J\wh{\bm{x}}_t\right),\bm{q}^T\bm{q}\right)_{L^2\left(\wh{\Omega}_h\right)}-\sum_{k}\frac{1}{2}\left\langle \bm{q}^*,A_n\bm{q}\right\rangle_{L^2(\p D^k)}.
\label{eq:rhs}
\end{align}
By the metric identity (\ref{eq:metric-id}), we have
\begin{align*}
\frac{\p}{\p\xi_1}A^1+\frac{\p}{\p\xi_2}A^2=\renewcommand\arraystretch{0.7}\begin{bmatrix}
\widehat{\nabla}\cdot\left(J\wh{\bm{x}}_t\right)&&\\
&\widehat{\nabla}\cdot\left(J\wh{\bm{x}}_t\right)&\\
&&\widehat{\nabla}\cdot\left(J\wh{\bm{x}}_t\right)
\end{bmatrix},
\end{align*}
such that the last volume term equals to
\begin{align*}
\frac{1}{2}\left(\wh{\nabla}\cdot\left(J\wh{\bm{x}}_t\right),\bm{q}^T\bm{q}\right)_{L^2\left(\wh{\Omega}_h\right)}=\frac{1}{2}\left(\left(\frac{\p}{\p\xi_1}A^1\right)\bm{q},\bm{q}\right)_{L^2\left(\wh{\Omega}_h\right)}
+\frac{1}{2}\left(\left(\frac{\p}{\p\xi_2}A^2\right)\bm{q},\bm{q}\right)_{L^2\left(\wh{\Omega}_h\right)}.
\end{align*}
Plugging this into (\ref{eq:rhs}), we see that all volume integral contributions vanish.

Now we plug in expression from (\ref{eq:ale-flux}) for numerical flux $\bm{q}^*$  and consider surface integral contributions
\begin{align}
RHS=-\sum_{k}\frac{1}{2}\left\langle \bm{q}^+-\tau_q A_n\jump{\bm{q}},A_n\bm{q}\right\rangle_{L^2(\p D^k)}.
\label{eq:surface-rhs}
\end{align}
By considering surface contributions from both neighboring elements, we obtain the desired theorem.
\begin{thm}
	The skew-symmetric ALE-DG formulation (\ref{eq:acous-ale-dg13}) using the DG method is energy stable in the following sense	
	\begin{equation}
	\begin{split}
	\frac{1}{2}\frac{\p}{\p\tau}\left(||p||_{J}^2+||u||_{J}^2+||v||_{J}^2\right)=-\tau_q[\![\bm{q}]\!]^TA_n^TA_n[\![\bm{q}]\!]\leq 0.
	\end{split}
	\end{equation}	
\end{thm}

To show near energy stability using weight-adjusted DG methods, we take
$$\bm{w} = \bm{q},\quad h = \frac{1}{2}\Pi_N\left(\widetilde{\bm{q}}^T\widetilde{\bm{q}}\right),$$
where
$$\widetilde{\bm{q}}=\renewcommand\arraystretch{0.7}\begin{bmatrix}
\widetilde{p}\\
\widetilde{u}\\
\widetilde{v}
\end{bmatrix},\quad \widetilde{p}=\frac{pJ}{J},\quad \widetilde{u}=\frac{uJ}{J},\quad \widetilde{v}=\frac{vJ}{J}.$$
In this case, (\ref{eq:rhs}) becomes
\begin{align}
RHS=-\frac{1}{2}\left(\wh{\nabla}\cdot\left(J\wh{\bm{x}}_t\right),\bm{q}^T\bm{q}\right)_{L^2\left(\wh{\Omega}_h\right)}+\frac{1}{2}\left(\wh{\nabla}\cdot\left(J\wh{\bm{x}}_t\right),\widetilde{\bm{q}}^T\widetilde{\bm{q}}\right)_{L^2\left(\wh{\Omega}_h\right)}-\sum_{k}\frac{1}{2}\left\langle \bm{q}^*,A_n\bm{q}\right\rangle_{L^2(\p D^k)}.
\label{eq:rhs-wadg}
\end{align}
Using the same approach in Section \ref{sec:energy-wadg}, we can bound (\ref{eq:rhs-wadg}) by
$$\Big|\Big|\frac{1}{2}\left(\wh{\nabla}\cdot\left(J\wh{\bm{x}}_t\right),\bm{q}^T\bm{q}\right)_{L^2(\Omega_h)}-\frac{1}{2}\left(\wh{\nabla}\cdot\left(J\wh{\bm{x}}_t\right),\widetilde{\bm{q}}^T\widetilde{\bm{q}}\right)_{L^2(\Omega_h)}\Big|\Big|\leq C(t) h^{2N+2}\leq C_{\text{max}}h^{2N+2}.$$
Combining this upper bound with surface integral contributions (\ref{eq:surface-rhs}), we obtain the following result.
\begin{thm}
	The skew-symmetric ALE-DG formulation (\ref{eq:acous-ale-dg13}) using the WADG method is energy stable up to a term which \textcolor{black}{converges} to zero under mesh refinement in the following sense	
	\begin{equation}
	\begin{split}
	\frac{1}{2}\frac{\p}{\p\tau}\left(||pJ||_{1/J}^2+||uJ||_{1/J}^2+||vJ||_{1/J}^2\right)\leq C h^{2N+2}-\tau_q[\![\bm{q}]\!]^TA_n^TA_n[\![\bm{q}]\!].
	\end{split}
	\end{equation}	
	\label{thm:wave-wadg}
\end{thm}
\section{Numerical experiments}\label{sec:ale-experiment}
In this section, we examine the accuracy and performance of the proposed method on moving meshes. We first run numerical experiments to verify the asymptotic energy conservation property discussed in Section \ref{sec:mesh-ale}. Then, we discuss the high order accuracy of our method for analytic solutions of the acoustic wave equation. Numerical experiments will be implemented on triangular meshes using degree $N$ polynomials.

\subsection{Energy conservation for mesh motion} 
In Section \ref{sec:mesh-ale}, we obtained the energy estimate (\ref{eq:estimate})
\begin{equation*}
\textcolor{black}{\Big|||uJ\left(\cdot,T\right)||_{1/J}^2-||uJ\left(\cdot,0\right)||_{1/J}^2\Big|\leq Ch^{2N+2}},
\end{equation*}
which indicates that, for any fixed $T$, the variation in energy will be bounded by some term that \textcolor{black}{converges to zero with the same rate as the optimal $L^2$ error estimate} as $h$ decreases. We  numerically verify this result in 1D and 2D using degree $N$ polynomial approximation spaces. 

\subsubsection{One-dimensional mesh motion}\label{sec:1d-mesh}
For our one-dimensional experiments, we choose the final time $T=0.5$ and the computational domain $[-1,1]$. We prescribe mesh motion by
\begin{align}
x=\xi+\frac{1}{4}\sin\left(\pi \tau\right)\left(1-\xi\right)\left(1+\xi\right).
\label{eq:1d-mesh-motion}
\end{align}
Under this mapping, the physical domain remains $[-1,1]$, but the interior coordinates change in time in a non-affine manner. We consider $\p u/\p t=0$ and the solution 
$$u(x,t)=\sin\left(\pi x\right),$$
with periodic boundary conditions. We compute the difference in energy between $t=0$ and $t=T$ for both standard DG methods (build and invert weighted mass matrices at each time step) and weight-adjusted DG (WADG) methods. We fix the order of approximation $N=2$ for the standard DG method. In all figures, the mesh size $h$ refers to the mesh size on the fixed reference domain $\wh{\Omega}$.

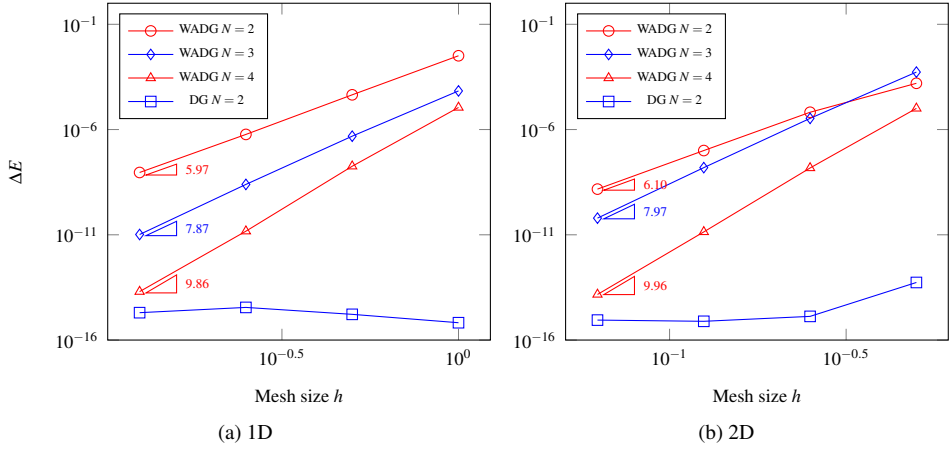
\begin{figure}[h]
	\centering
	\hspace{-4em}
	\subfloat[1D]{
		\begin{tikzpicture}
		\begin{loglogaxis}[
		width=0.55\textwidth,
		height=0.5\textwidth,
		title style = {font=\large},
		xlabel= Mesh size $h$,
		ylabel= $\Delta E$,
		label style = {font=\scriptsize},
		ticklabel style = {font=\scriptsize},
		ymin=1e-16,
		ymax=1e0,
		legend style={font=\tiny},
		legend pos = north west
		]
		%T=0.5
		\addplot[color=red,mark=o] coordinates {
			(1, 3.1995e-03)
			(0.5, 4.4645e-05)
			(0.25, 5.8666e-07)
			(0.125,9.1273e-09)
		};
		
		\addplot[color=blue,mark=diamond] coordinates {
			(1,   6.8349e-05)
			(0.5, 4.8430e-07)
			(0.25, 2.4742e-09 )
			(0.125,1.0331e-11)
		};
		\addplot[color=red,mark=triangle] coordinates {
			(1, 1.1190e-05)
			(0.5, 1.8051e-08)
			(0.25, 1.4819e-11 )
			(0.125,1.9984e-14)
		};
		
		\addplot[color=blue,mark=square] coordinates {
			(1, 6.6613e-16)
			(0.5, 1.6653e-15)
			(0.25, 3.5527e-15)
			(0.125,1.9984e-15)
		};
		
		\logLogSlopeTriangle{0.18}{0.08}{0.49}{5.97}{red};
		\logLogSlopeTriangle{0.18}{0.08}{0.31}{7.87}{blue};
		\logLogSlopeTriangle{0.18}{0.08}{0.14}{9.86}{red};
		
		\legend{WADG $N=2$,WADG $N=3$,WADG $N=4$,DG $N=2$}
		\end{loglogaxis}
		\end{tikzpicture}
		\label{fig:energy-1d}	
	}
	\subfloat[2D]{
		\begin{tikzpicture}
		\begin{loglogaxis}[
		width=0.55\textwidth,
		height=0.5\textwidth,
		title style = {font=\normalsize},
		xlabel= Mesh size $h$,
		label style = {font=\scriptsize},
		ticklabel style = {font=\scriptsize},
		ymin=1e-16,
		ymax=1e0,
		legend style={font=\tiny},
		legend pos = north west
		]
		%T=0.5 CFL=0.1
		\addplot[color=red,mark=o] coordinates {
			(0.5, 1.5749e-04)
			(0.25, 6.6328e-06 )
			(0.125, 1.0063e-07)
			(0.0625,1.4796e-09)
		};
		
		\addplot[color=blue,mark=diamond] coordinates {
			(0.5, 5.3681e-04)
			(0.25, 3.3963e-06 )
			(0.125, 1.5505e-08)
			(0.0625,6.1953e-11)
		};
		\addplot[color=red,mark=triangle] coordinates {
			(0.5, 1.0125e-05)
			(0.25,  1.4948e-08)
			(0.125, 1.3717e-11)
			(0.0625,1.4655e-14)
		};
		
		\addplot[color=blue,mark=square] coordinates {
			(0.5, 5.4401e-14)
			(0.25, 1.3323e-15)
			(0.125, 7.7716e-16)
			(0.0625,8.8818e-16)
		};
		
		\logLogSlopeTriangle{0.18}{0.08}{0.445}{6.10}{red};
		\logLogSlopeTriangle{0.18}{0.08}{0.36}{7.97}{blue};
		\logLogSlopeTriangle{0.18}{0.08}{0.135}{9.96}{red};
		\legend{WADG $N=2$,WADG $N=3$,WADG $N=4$,DG $N=2$}
		\end{loglogaxis}
		\end{tikzpicture}
		\label{fig:energy-2d-curved}
	}
	\caption[Energy variation of sufficiently regular solutions for different orders of approximation.]{Energy variation of sufficiently regular solutions for different orders of approximation.}
\end{figure}

From Figure \ref{fig:energy-1d}, we observe that the convergence behavior of energy variation $\Delta E$ is exactly what we expect from Theorem~\ref{thm:upper-bound}. For WADG, we achieve the convergence rate $O(h^{2N+2})$ predicted in (\ref{eq:estimate}). We also include the change in energy for the standard DG method to demonstrate that the skew-symmetric semi-discretely formulation (\ref{eq:skew-mesh}) is energy conservative if we use the standard weighted mass matrices. Note that the change in energy is not exactly zero because the proof of energy conservation is semi-discrete, not fully discrete. The energy variation can be further reduced if we decrease the time step size, or if relaxation Runge-Kutta methods are used \cite{ketcheson2019relaxation}. Thus, for the one-dimensional case, our numerical results are consistent with Theorem \ref{thm:energy-conservation} and \ref{thm:upper-bound} on energy conservation.

\subsubsection{Two-dimensional mesh motion}\label{sec:2d-mesh}
In our 2D experiments, we choose the final time $T=0.5$ and the computational domain $[-1,1]\times[-1,1]$. We run experiments using the following mesh motion
\begin{equation}
\begin{split}
x_1&=\xi_1 + \frac{1}{4}\sin\left(\pi \tau\right)\sin\left(\pi \xi_1\right)\left(1-\xi_1\right)\left(1+\xi_1\right),\\ x_2&=\xi_2 + \frac{1}{4}\sin\left(\pi \tau\right)\sin\left(\pi \xi_2 \right)\left(1-\xi_2\right)\left(1+\xi_2\right),
\end{split}
\label{eq:2d-mesh-motion}
\end{equation}
which results in time-dependent curvilinear meshes. Again, the physical domains remain the same as the reference domain. We consider $\p u/\p t=0$ and the constant solution
$$u(x_1,x_2,t)=\sin\left(\pi x_1\right)\cos\left(\pi x_2\right),$$
with periodic boundary conditions.

From Figure \ref{fig:energy-2d-curved}, we can see that the convergence rate of the energy variation is consistent with the energy estimate (\ref{eq:estimate}). The energy variation for standard DG method with $N=2$ is presented as well, and the result is smaller than $O(10^{-12})$, which is very close to machine epsilon. 

\subsection{ALE-DG for wave propagation}
In this section, we investigate the accuracy of the ALE-DG method for analytic solutions of the acoustic wave equation. We implemented two-dimensional simulations of the acoustic wave equation on moving meshes to check convergence behavior and energy stability. We simulate wave propagation until final time $T=1.5$ with mesh motion given by (\ref{eq:2d-mesh-motion}).

%\subsubsection{Simulations in homogeneous media}

In this setting, we assume homogeneous media such that $c=1$. We choose the pressure $p\left(\bm{x},t\right)$ to be
$$p\left(\bm{x},t\right)=\sin\left(\pi x_1\right)\sin\left(\pi x_2\right)\cos\left(\sqrt{2}\pi t\right),$$
and solve for $u\left(\bm{x},t\right),v\left(\bm{x},t\right)$ which satisfy (\ref{eq:conservative})
\begin{align*}
u\left(\bm{x},t\right)&=-\frac{\sqrt{2}}{2}\cos\left(\pi x_1\right)\sin\left(\pi x_2\right)\sin\left(\sqrt{2}\pi t\right)\\
v\left(\bm{x},t\right)&=-\frac{\sqrt{2}}{2}\sin\left(\pi x_1\right)\cos\left(\pi x_2\right)\sin\left(\sqrt{2}\pi t\right)
\end{align*}
\begin{figure}[h]
	\centering
	\hspace{-4em}
	\subfloat[Central flux ($\tau_q=0$)]{
		\begin{tikzpicture}
		\begin{loglogaxis}[
		width=0.55\textwidth,
		height=0.5\textwidth,
		title style = {font=\large},
		xlabel= Mesh size $h$,
		ylabel= $|\Delta E|$,
		label style= {font=\scriptsize},
		ticklabel style = {font=\scriptsize},
		ymax = 0.5e-00,
		ymin = 5e-16,
		legend style={font=\tiny},
		legend pos = north west
		]
		
		%T=1.5
		\addplot[color=red,mark=o] coordinates {
			(0.5,  1.0360e-03)
			(0.25, 1.6467e-05)
			(0.125,1.6340e-07)
			(0.0625,2.6095e-09)
		};
		
		\addplot[color=blue,mark=diamond] coordinates {
			(0.5, 8.5672e-04)
			(0.25, 3.5489e-06)
			(0.125,1.8629e-08)
			(0.0625,7.7776e-11)
		};
		\addplot[color=red,mark=triangle] coordinates {
			(0.5,  1.1020e-05)
			(0.25, 2.0094e-08 )
			(0.125,2.0364e-11 )
			(0.0625,2.2871e-14)
		};
		
		%DG N=2
		\addplot[color=blue,mark=square] coordinates {
			(0.5, 8.4377e-15)
			(0.25,  9.4369e-15)
			(0.125,1.9984e-15)
			(0.0625,4.3299e-15)
		};

		\logLogSlopeTriangle{0.18}{0.08}{0.445}{6.09}{red};
		\logLogSlopeTriangle{0.18}{0.08}{0.345}{7.98}{blue};
		\logLogSlopeTriangle{0.18}{0.08}{0.11}{10.21}{red};
		\legend{WADG $N=2$,WADG $N=3$,WADG $N=4$, DG $N=2$}
		\end{loglogaxis}
		\end{tikzpicture}
		\label{fig:energy-central-homo}	
	}
	%\hspace{0.2cm}
	\subfloat[Penalty flux ($\tau_q=1$)]{
		\begin{tikzpicture}
		\begin{loglogaxis}[
		width=0.55\textwidth,
		height=0.5\textwidth,
		title style = {font=\large},
		xlabel= Mesh size $h$,
		label style = {font=\scriptsize},
		ticklabel style = {font=\scriptsize},
		ymax = 0.5e-00,
		ymin = 5e-16,
		legend style={font=\tiny},
		legend pos = south east
		]
	
		%T=1.5
		\addplot[color=red,mark=o] coordinates {
			(0.5, 7.0581e-02)
			(0.25, 4.8919e-03 )
			(0.125, 1.8978e-04)
			(0.0625, 7.4874e-06)
		};
		
		\addplot[color=blue,mark=diamond] coordinates {
			(0.5,8.1823e-03)
			(0.25, 1.1266e-04 )
			(0.125,2.0229e-06)
			(0.0625,3.4381e-08)
		};
		\addplot[color=red,mark=triangle] coordinates {
			(0.5, 7.1003e-04)
			(0.25, 4.1058e-06 )
			(0.125,1.0982e-08 )
			(0.0625,3.0910e-11)
		};
		%DG N=2
		\addplot[color=blue,mark=square] coordinates {
			(0.5, 7.1391e-02)
			(0.25,  4.9126e-03)
			(0.125,1.9009e-04)
			(0.0625,7.4920e-06)
		};
		
		\logLogSlopeTriangle{0.18}{0.08}{0.67}{4.66}{red};
		\logLogSlopeTriangle{0.18}{0.08}{0.52}{5.87}{blue};
		\logLogSlopeTriangle{0.18}{0.08}{0.32}{8.47}{red};
		\legend{WADG $N=2$,WADG $N=3$,WADG $N=4$, DG $N=2$}
		\end{loglogaxis}
		\end{tikzpicture}
		\label{fig:energy-penalty-homo}
	}
	\caption[Energy variation of DG and WADG for different orders of approximation]{Energy variation of DG and WADG for different orders of approximation}
\end{figure}
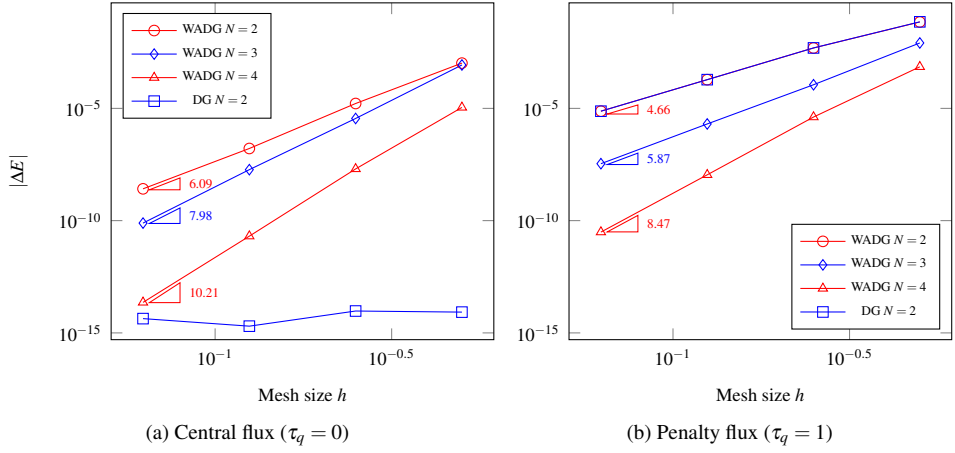

In the first experiment, we check energy stability for both central fluxes and penalty fluxes by setting the penalty parameter to $\tau_q=0$ and $\tau_q=1$. From Figure \ref{fig:energy-central-homo}, we observe the expected $O(h^{2N+2})$ convergence rate for $\Delta E$. In Figure \ref{fig:energy-penalty-homo}, the result is slightly different from Figure \ref{fig:energy-central-homo}, as we add energy dissipation in this case. We can see that the effect of energy dissipation becomes smaller and smaller as we decrease the mesh size. This makes sense because the dissipation is proportional to solution jumps across element interfaces, which decreases as the mesh size decreases. \textcolor{black}{Moreover, we can find that the dissipation dominates the change in energy due to WADG methods.} Similar to previous results for solutions of $\p u/\p t=0$ under mesh motion, the proposed skew-symmetric WADG method is again energy conservative for the wave equation up to a term which \textcolor{black}{converges to zero with the same rate as the optimal $L^2$ error estimate} under mesh refinement.

\begin{figure}[h]
	\centering
	\hspace{-4em}
	\subfloat[Central flux ($\tau_q=0$)]{
		\begin{tikzpicture}
		\begin{loglogaxis}[
		width=0.55\textwidth,
		height=0.5\textwidth,
		title style = {font=\large},
		xlabel= Mesh size $h$,
		ylabel= $L^2$ error,
		label style = {font=\scriptsize},
		ticklabel style = {font=\scriptsize},
		ymax = 1e-00,
		ymin = 0.5e-7,
		legend style={font=\tiny},
		legend pos = south east
		]
		
		%T=1.5
		\addplot[color=red,mark=o] coordinates {
			(0.5, 1.2883e-01)
			(0.25, 2.2103e-02 )
			(0.125,2.9360e-03)
			(0.0625,3.7840e-04)
		};
		
		\addplot[color=blue,mark=diamond] coordinates {
			(0.5, 4.0139e-02)
			(0.25, 2.9747e-03 )
			(0.125,3.0142e-04 )
			(0.0625,2.2626e-05)
		};
		\addplot[color=red,mark=triangle] coordinates {
			(0.5, 8.4904e-03)
			(0.25, 4.5340e-04)
			(0.125,2.1746e-05 )
			(0.0625,8.4503e-07)
		};
		
		\logLogSlopeTriangle{0.18}{0.08}{0.53}{2.95}{red};
		\logLogSlopeTriangle{0.18}{0.08}{0.36}{3.75}{blue};
		\logLogSlopeTriangle{0.18}{0.08}{0.17}{4.70}{red};
		\legend{$N=2$,$N=3$,$N=4$}
		\end{loglogaxis}
		\end{tikzpicture}
		\label{fig:convergence-central-homo}	
	}
	\subfloat[Penalty flux ($\tau_q=1$)]{
		\begin{tikzpicture}
		\begin{loglogaxis}[
		width=0.55\textwidth,
		height=0.5\textwidth,
		title style = {font=\large},
		xlabel= Mesh size $h$,
		label style = {font=\scriptsize},
		ticklabel style = {font=\scriptsize},
		ymax = 1e-00,
		ymin = 0.5e-7,
		legend style={font=\tiny},
		legend pos = south east
		]

		%T=1.5
		\addplot[color=red,mark=o] coordinates {
			(0.5, 1.0595e-01 )
			(0.25, 1.0533e-02)
			(0.125,9.0159e-04)
			(0.0625,9.7360e-05)
		};
		
		\addplot[color=blue,mark=diamond] coordinates {
			(0.5, 2.8987e-02)
			(0.25, 1.7870e-03)
			(0.125,1.2602e-04)
			(0.0625,8.6660e-06)
		};
		\addplot[color=red,mark=triangle] coordinates {
			(0.5, 4.6943e-03)
			(0.25, 1.8442e-04 )
			(0.125,6.1508e-06)
			(0.0625,1.9523e-07)
		};
		
		\logLogSlopeTriangle{0.18}{0.08}{0.45}{3.21}{red};
		\logLogSlopeTriangle{0.18}{0.08}{0.305}{3.85}{blue};
		\logLogSlopeTriangle{0.18}{0.08}{0.08}{4.96}{red};
		\legend{$N=2$,$N=3$,$N=4$}
		\end{loglogaxis}
		\end{tikzpicture}
		\label{fig:convergence-penalty-homo}
	}
	\caption[Convergence behavior of WADG for different orders of approximation]{Convergence behavior of WADG for different orders of approximation\label{fig:convergence-homo}}
\end{figure}
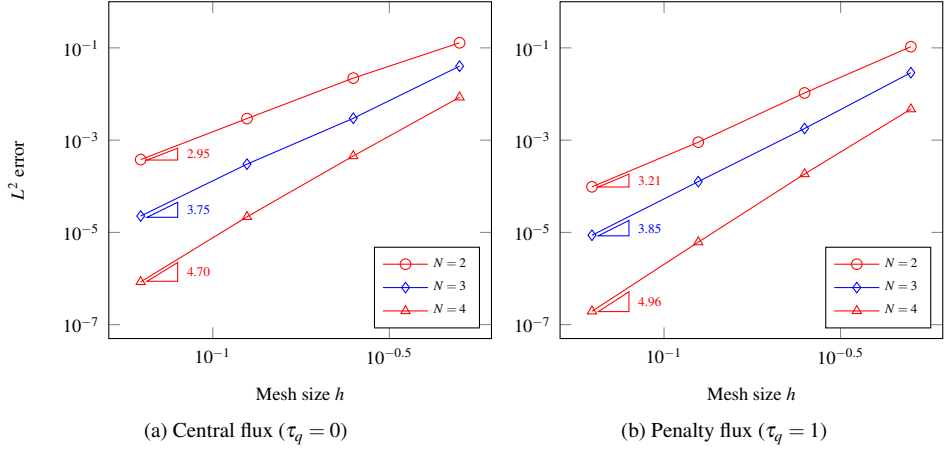
\begin{figure}[h]
	\centering
	\subfloat[$T=0.5$ (moving mesh)]{
		\includegraphics[width=0.45\linewidth,height=0.35\linewidth]{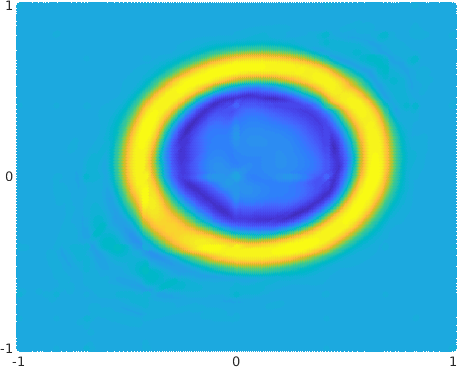}
	}
	\subfloat[$T=0.5$ (static mesh)]{
		\includegraphics[width=0.45\linewidth,height=0.35\linewidth]{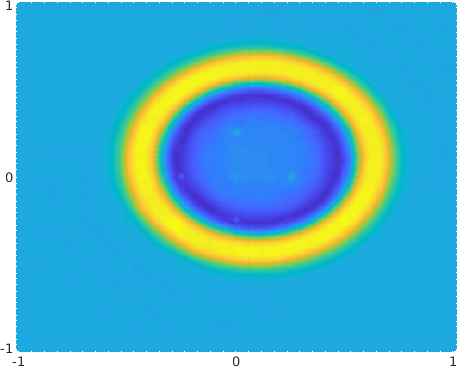}
	}
\\
	\subfloat[$T=0.9$ (moving mesh)]{
		\includegraphics[width=0.45\linewidth,height=0.35\linewidth]{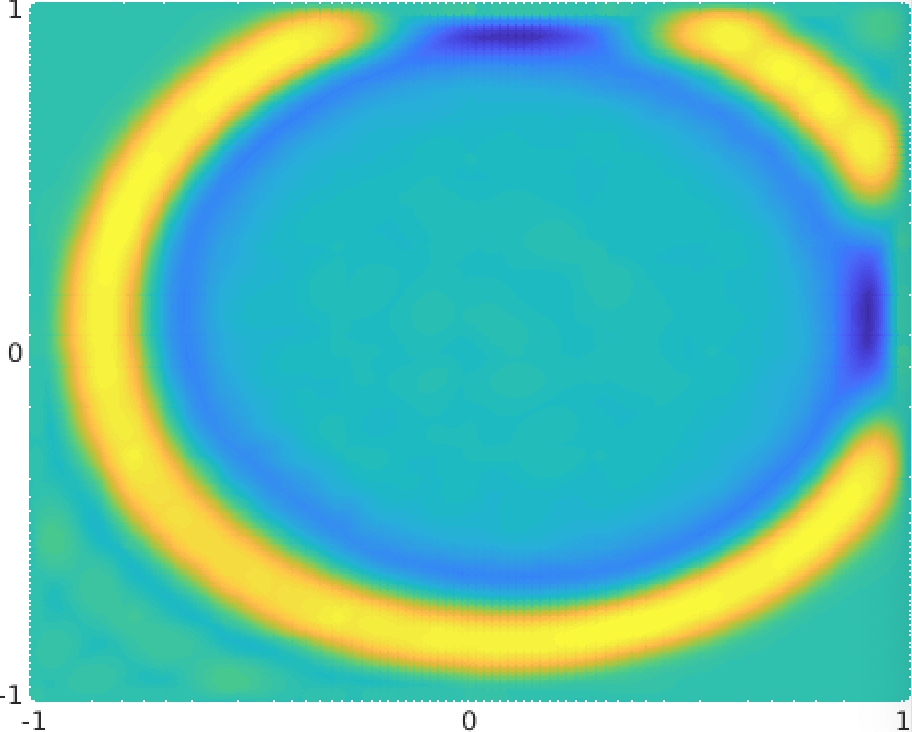}
	}
	\subfloat[$T=0.9$ (static mesh)]{
		\includegraphics[width=0.45\linewidth,height=0.35\linewidth]{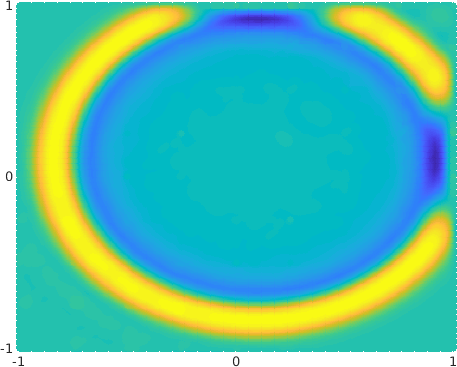}
	}		
	\caption[Wave propagation triggered by a Gaussian pulse]{Wave propagation triggered by a Gaussian pulse}
	\label{fig:wave-propagation}
\end{figure}

\begin{figure}[h]
	\centering
	\subfloat[Moving mesh at $T=0$]{
		\includegraphics[width=0.3\linewidth,height=0.25\linewidth]{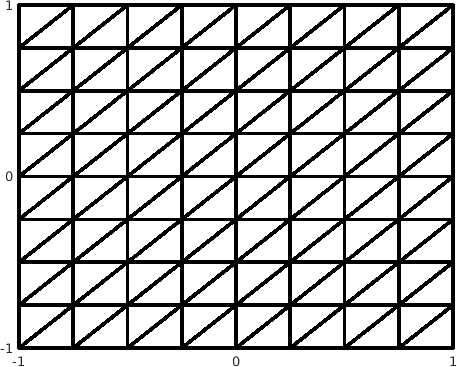}
	}
	\subfloat[Moving mesh at $T=0.5$]{
		\includegraphics[width=0.3\linewidth,height=0.25\linewidth]{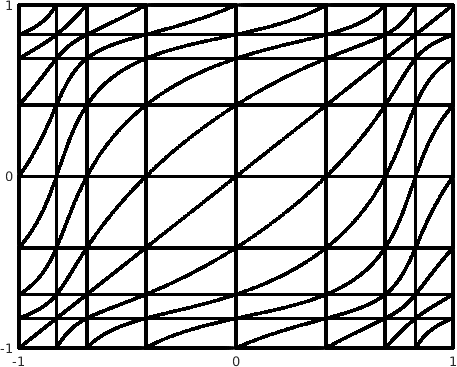}
	}
	\subfloat[Moving mesh at $T=0.9$]{
		\includegraphics[width=0.3\linewidth,height=0.25\linewidth]{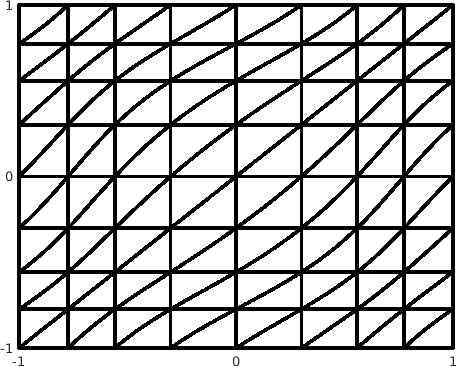}
	}
	\caption[Moving mesh at different time]{Moving mesh at different time}
	\label{fig:wave-propagation-mesh}
\end{figure}

In the second experiment, we check the convergence of the $L^2$ error for both central and penalty fluxes (see Figure~\ref{fig:convergence-homo}). The rate of convergence for penalty fluxes is close to the theoretical rate $O(h^{N+\frac{1}{2}})$ \cite{johnson1986analysis}, which is faster than that for central fluxes. To better illustrate performance of the proposed method, we simulate wave propagation triggered by a Gaussian pulse at point $\bm{x}=(0.1,0.1)$ as follows
$$p_0(x_1,x_2) = e^{-100\left((x-0.1)^2+(y-0.1)^2\right)}.$$
The wave simulation is conducted on a static domain with
a moving mesh, whose motion is given by (\ref{eq:2d-mesh-motion}). We will also simulate on the same domain with a static mesh for comparison.  Figure~\ref{fig:wave-propagation} shows the pressure $p$ at times $T=0.5$ and $T=0.9$ on both moving and static mesh. We observe that the numerical result on a moving mesh is slightly underresolved near wave fronts, and this can be eliminated as we increase the order of approximation. Overall, our proposed method can easily handle mesh motion and the yielding results agree with the reference solution on a static mesh.

\subsection{Energy conservation for discontinuous solutions}
In this section, we investigate energy conservation for discontinuous or under-resolved solutions on moving meshes. We first consider $\p u/\p t=0$ and constant solutions
$$u(x,t)=\begin{cases}
\sin\left(\pi x\right), &\mbox{if}\ x\leq 0,\\
\cos\left(\pi x\right), &\mbox{if}\ x> 0,
\end{cases}$$
for 1D experiment and
$$u(x_1,x_2,t)=\begin{cases}
\sin\left(\pi x_1\right)\cos\left(\pi x_2\right), &\mbox{if}\ x_1\leq 0,\\
\cos\left(\pi x_1\right)\cos\left(\pi x_2\right), &\mbox{if}\ x_1> 0,
\end{cases}$$
for 2D experiment. We choose the final time $T=0.5$ and periodic boundary conditions. We use the same computational domains and mesh motion in Section~\ref{sec:1d-mesh} and \ref{sec:2d-mesh}.

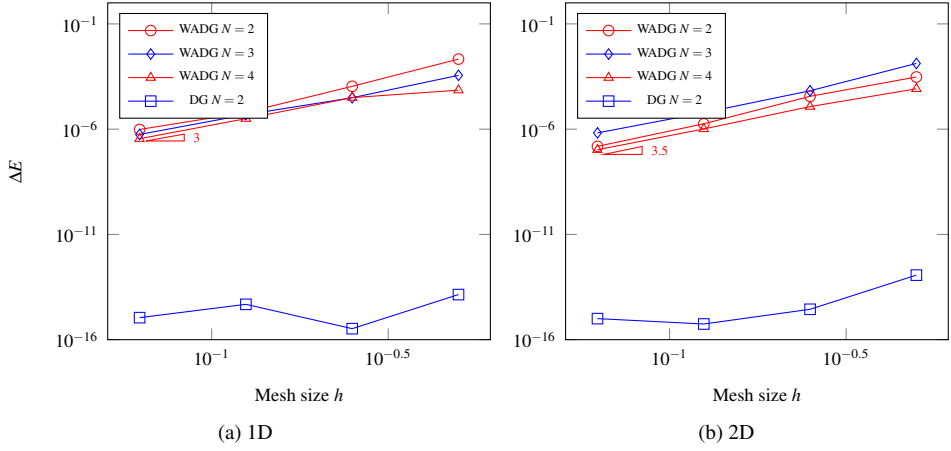
\begin{figure}[h]
	\centering
	\hspace{-4em}
	\subfloat[1D]{
		\begin{tikzpicture}
		\begin{loglogaxis}[
		width=0.55\textwidth,
		height=0.5\textwidth,
		title style = {font=\large},
		xlabel= Mesh size $h$,
		ylabel= $\Delta E$,
		label style = {font=\scriptsize},
		ticklabel style = {font=\scriptsize},
		ymin=1e-16,
		ymax=1e-0,
		legend style={font=\tiny},
		legend pos = north west
		]
		%T=0.5
		\addplot[color=red,mark=o] coordinates {
			(0.5, 2.1345e-03)
			(0.25,1.0875e-04 )
			(0.125, 6.2530e-06)
			(0.0625,9.6779e-07)
		};
		
		\addplot[color=blue,mark=diamond] coordinates {
			(0.5,  3.6048e-04)
			(0.25, 3.1112e-05)
			(0.125, 4.9018e-06)
			(0.0625,5.5990e-07)
		};
		\addplot[color=red,mark=triangle] coordinates {
			(0.5, 7.1299e-05)
			(0.25, 3.1479e-05)
			(0.125, 3.1554e-06)
			(0.0625,3.5184e-07)
		};
		
		\addplot[color=blue,mark=square] coordinates {
			(0.5, 1.3878e-14)
			(0.25,  3.3307e-16 )
			(0.125, 4.7740e-15)
			(0.0625,1.1102e-15)
		};

		\logLogSlopeTriangle{0.2}{0.1}{0.59}{3}{red};
		
		\legend{WADG $N=2$,WADG $N=3$,WADG $N=4$,DG $N=2$}
		\end{loglogaxis}
		\end{tikzpicture}
		%\label{fig:energy-1d-dis}	
	}
	\subfloat[2D]{
		\begin{tikzpicture}
		\begin{loglogaxis}[
		width=0.55\textwidth,
		height=0.5\textwidth,
		title style = {font=\normalsize},
		xlabel= Mesh size $h$,
		%ylabel= $L^2$ error,
		label style = {font=\scriptsize},
		ticklabel style = {font=\scriptsize},
		ymin=1e-16,
		ymax=1e-0,
		legend style={font=\tiny},
		legend pos = north west
		]
	
		\addplot[color=red,mark=o] coordinates {
			(0.5, 3.0099e-04 )
			(0.25, 3.6313e-05)
			(0.125, 1.7988e-06)
			(0.0625,1.5250e-07)
		};
		
		\addplot[color=blue,mark=diamond] coordinates {
			(0.5, 1.3262e-03)
			(0.25, 6.6716e-05)
			(0.125, 6.0515e-06)
			(0.0625,6.6526e-07)
		};
		\addplot[color=red,mark=triangle] coordinates {
			(0.5, 8.2976e-05)
			(0.25,  1.1815e-05)
			(0.125, 1.0240e-06)
			(0.0625,1.0528e-07)
		};
		
		\addplot[color=blue,mark=square] coordinates {
			(0.5, 1.1657e-13)
			(0.25, 2.7756e-15)
			(0.125, 5.5511e-16 )
			(0.0625,9.9920e-16)
		};
		\logLogSlopeTriangle{0.2}{0.1}{0.55}{3.5}{red};
		
		\legend{WADG $N=2$,WADG $N=3$,WADG $N=4$,DG $N=2$}
		\end{loglogaxis}
		\end{tikzpicture}
		%\label{fig:energy-2d-dis-curved}
	}
	
	\caption[Energy variation of discontinuous constant solutions for different orders of approximation.]{Energy variation of discontinuous constant solutions for different orders of approximation.}
	\label{fig:energy-dis-mesh}
\end{figure}

In Section~\ref{sec:energy-wadg}, for sufficiently regular solution $u$, we have shown that there exists an upper bound for the energy variation using WADG method given by
$$\textcolor{black}{\Big|||uJ\left(\cdot,T\right)||_{1/J}^2-||uJ\left(\cdot,0\right)||_{1/J}^2\Big|\leq Ch^{r}},$$
where $r=2N+2$. However, from Figure~\ref{fig:energy-dis-mesh}, when solutions are not sufficient regular, the rate $r$ is much lower than that for sufficiently regular solutions and becomes independent of the order of approximation $N$. When we use the standard DG method,  the skew-symmetric semi-discretely formulation (\ref{eq:skew-mesh}) remains energy conservative for discontinuous solutions.

\begin{figure}[h]
	\centering
	\hspace{-4em}
	\subfloat[1D]{
		\begin{tikzpicture}
		\begin{loglogaxis}[
		width=0.55\textwidth,
		height=0.5\textwidth,
		title style = {font=\large},
		xlabel= Mesh size $h$,
		ylabel= $\Delta E$,
		label style = {font=\scriptsize},
		ticklabel style = {font=\scriptsize},
		ymin=1e-12,
		ymax=1e0,
		legend style={font=\tiny},
		legend pos = north west
		]
		%T=1.5 CFL=0.0001
		\addplot[color=red,mark=o] coordinates {
			(0.5, 7.5733e-04)
			(0.25,  4.1125e-05 )
			(0.125, 4.8094e-06 )
			(0.0625,7.7929e-07)
		};
		
		\addplot[color=blue,mark=diamond] coordinates {
			(0.5, 1.6188e-04)
			(0.25, 2.0460e-05)
			(0.125, 3.4418e-06)
			(0.0625,3.4220e-07)
		};
		\addplot[color=red,mark=triangle] coordinates {
			(0.5, 1.2883e-04)
			(0.25, 1.5389e-05)
			(0.125, 1.2681e-06)
			(0.0625,1.7046e-07)
		};
		
		\addplot[color=blue,mark=square] coordinates {
			(0.5, 3.3568e-11)
			(0.25,  5.4636e-12)
			(0.125, 1.1493e-11)
			(0.0625,1.1855e-11)
		};
		
		\logLogSlopeTriangle{0.2}{0.1}{0.43}{3}{red};
		
		\legend{WADG $N=2$,WADG $N=3$,WADG $N=4$,DG $N=2$}
		\end{loglogaxis}
		\end{tikzpicture}
		\label{fig:energy-1d-dis-wave}	
	}
	\subfloat[2D]{
		\begin{tikzpicture}
		\begin{loglogaxis}[
		width=0.55\textwidth,
		height=0.5\textwidth,
		title style = {font=\normalsize},
		xlabel= Mesh size $h$,
		label style = {font=\scriptsize},
		ticklabel style = {font=\scriptsize},
		ymin=1e-12,
		ymax=1e0,
		legend style={font=\tiny},
		legend pos = north west
		]
		%T=1.5 CFL=0.1
		\addplot[color=red,mark=o] coordinates {
			(0.5, 8.1607e-03)
			(0.25, 4.8780e-04)
			(0.125, 4.8831e-05)
			(0.0625,4.5455e-06)
		};
		
		\addplot[color=blue,mark=diamond] coordinates {
			(0.5, 1.8592e-02)
			(0.25, 1.3885e-03)
			(0.125, 1.6574e-04)
			(0.0625,1.7793e-05)
		};
		\addplot[color=red,mark=triangle] coordinates {
			(0.5, 9.3607e-04)
			(0.25,  9.1200e-05)
			(0.125, 8.7387e-06)
			(0.0625,7.1626e-07)
		};
		%CFL=0.01
		\addplot[color=blue,mark=square] coordinates {
			(0.5, 1.8651e-10)
			(0.25, 5.7757e-11)
			(0.125, 3.4651e-11)
			(0.0625,2.8018e-11)
		};
		
		\logLogSlopeTriangle{0.2}{0.1}{0.48}{3.5}{red};
		\legend{WADG $N=2$,WADG $N=3$,WADG $N=4$,DG $N=2$}
		\end{loglogaxis}
		\end{tikzpicture}
		\label{fig:energy-2d-dis-wave}
	}
	\caption[Energy variation of discontinuous initial conditions for different orders of approximation (central fluxes).]{Energy variation of discontinuous initial conditions for different orders of approximation (central fluxes).}
	\label{fig:energy-dis-wave}
\end{figure}
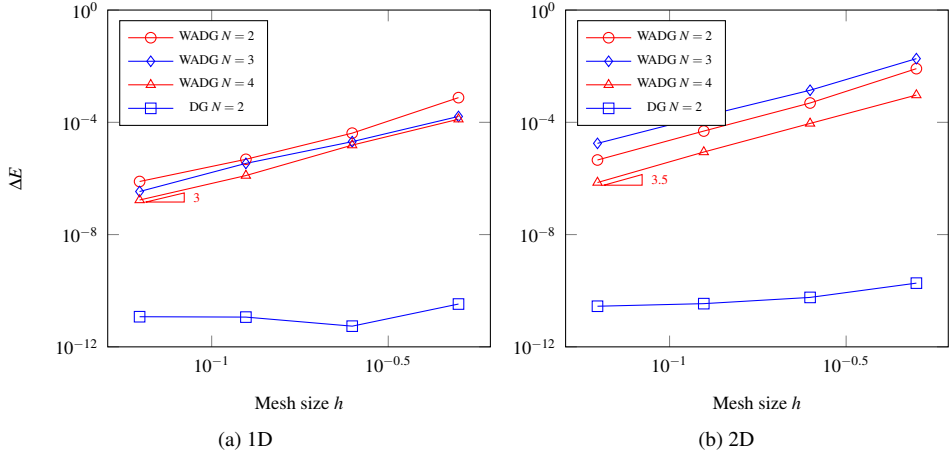

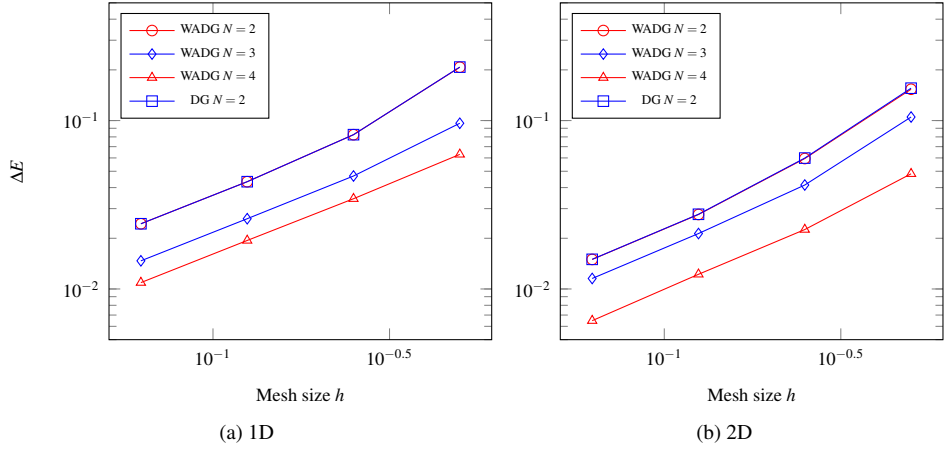
\begin{figure}[h]
	\centering
	\hspace{-4em}
	\subfloat[1D]{
		\begin{tikzpicture}
		\begin{loglogaxis}[
		width=0.55\textwidth,
		height=0.5\textwidth,
		title style = {font=\large},
		xlabel= Mesh size $h$,
		ylabel= $\Delta E$,
		label style = {font=\scriptsize},
		ticklabel style = {font=\scriptsize},
		ymin=5e-03,
		ymax=5e-01,
		legend style={font=\tiny},
		legend pos = north west
		]
		%T=1.5 CFL=0.01
		\addplot[color=red,mark=o] coordinates {
			(0.5, 2.0740e-01)
			(0.25,  8.2438e-02)
			(0.125, 4.3351e-02)
			(0.0625,2.4390e-02)
		};
		
		\addplot[color=blue,mark=diamond] coordinates {
			(0.5, 9.6406e-02 )
			(0.25, 4.6862e-02 )
			(0.125,2.6131e-02 )
			(0.0625,1.4694e-02)
		};
		\addplot[color=red,mark=triangle] coordinates {
			(0.5,  6.2982e-02 )
			(0.25, 3.4350e-02)
			(0.125, 1.9414e-02)
			(0.0625,1.0911e-02)
		};
		
		\addplot[color=blue,mark=square] coordinates {
			(0.5, 2.0777e-01)
			(0.25,  8.2460e-02)
			(0.125, 4.3354e-02)
			(0.0625,2.4390e-02)
		};

		\legend{WADG $N=2$,WADG $N=3$,WADG $N=4$,DG $N=2$}
		\end{loglogaxis}
		\end{tikzpicture}
		\label{fig:energy-1d-dis-wave-penalty}	
	}
	\subfloat[2D]{
		\begin{tikzpicture}
		\begin{loglogaxis}[
		width=0.55\textwidth,
		height=0.5\textwidth,
		title style = {font=\normalsize},
		xlabel= Mesh size $h$,
		%ylabel= $L^2$ error,
		label style = {font=\scriptsize},
		ticklabel style = {font=\scriptsize},
		ymin=5e-03,
		ymax=5e-01,
		legend style={font=\tiny},
		legend pos = north west
		]
		%T=1.5 CFL=0.1
		\addplot[color=red,mark=o] coordinates {
			(0.5, 1.5347e-01)
			(0.25, 5.9501e-02)
			(0.125, 2.7686e-02)
			(0.0625,1.4976e-02)
		};
		
		\addplot[color=blue,mark=diamond] coordinates {
			(0.5, 1.0497e-01)
			(0.25, 4.1430e-02)
			(0.125, 2.1338e-02)
			(0.0625,1.1542e-02)
		};
		\addplot[color=red,mark=triangle] coordinates {
			(0.5, 4.8311e-02)
			(0.25,  2.2534e-02)
			(0.125, 1.2222e-02)
			(0.0625,6.4832e-03)
		};
		%CFL=0.01
		\addplot[color=blue,mark=square] coordinates {
			(0.5, 1.5571e-01)
			(0.25,  5.9878e-02)
			(0.125, 2.7773e-02)
			(0.0625,1.4997e-02)
		};
		
		\legend{WADG $N=2$,WADG $N=3$,WADG $N=4$,DG $N=2$}
		\end{loglogaxis}
		\end{tikzpicture}
		\label{fig:energy-2d-dis-wave-penalty}
	}
	\caption[Energy variation of discontinuous initial conditions for different orders of approximation (penalty fluxes).]{Energy variation of discontinuous initial conditions for different orders of approximation (penalty fluxes).}
	\label{fig:energy-dis-wave-penalty}
\end{figure}
Next, we simulate 1D and 2D wave propagation with discontinuous initial conditions and reflective boundary conditions. We set initial conditions in 1D to be
\begin{align*}
p_0=\begin{cases*}
\sin\left(\pi x\right),\quad \text{if}\ x\leq 0,\\
\cos\left(\pi x\right),\quad \text{if}\ x> 0,
\end{cases*},\qquad u_0 = 0,
\end{align*}
and in 2D to be
\begin{align*}
p_0=\begin{cases*}
\sin\left(\pi x_1\right)\cos\left(\pi x_2\right),\quad \text{if}\ x\leq 0,\\
\cos\left(\pi x_1\right)\cos\left(\pi x_2\right),\quad \text{if}\ x> 0,
\end{cases*},\qquad u_0 = 0,\qquad v_0 = 0.
\end{align*}
\begin{figure}[h]
	\centering
	\subfloat[DG (central flux)]{
		\includegraphics[width=0.45\linewidth,height=0.3\linewidth]{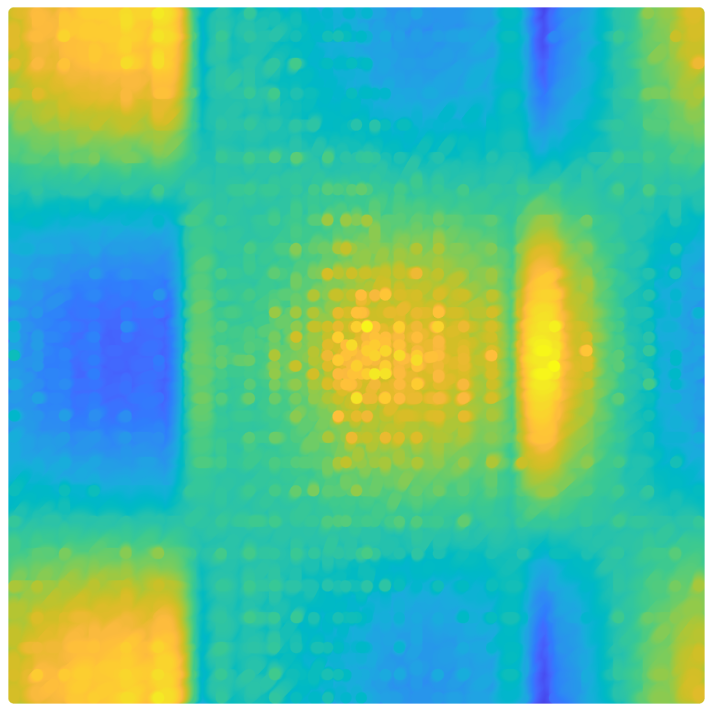}
	}
	\subfloat[DG (penalty flux)]{
		\includegraphics[width=0.45\linewidth,height=0.3\linewidth]{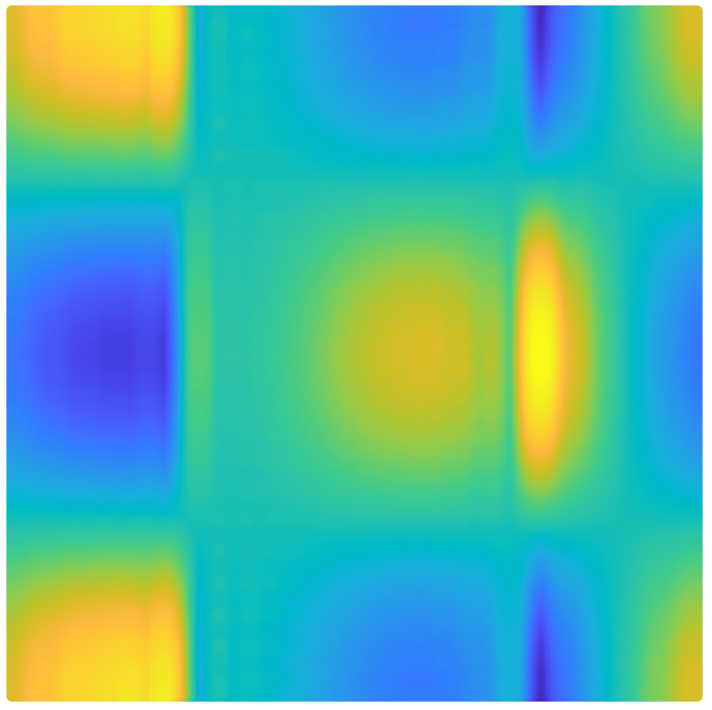}
	}\\
	\subfloat[WADG (central flux)]{
		\includegraphics[width=0.45\linewidth,height=0.3\linewidth]{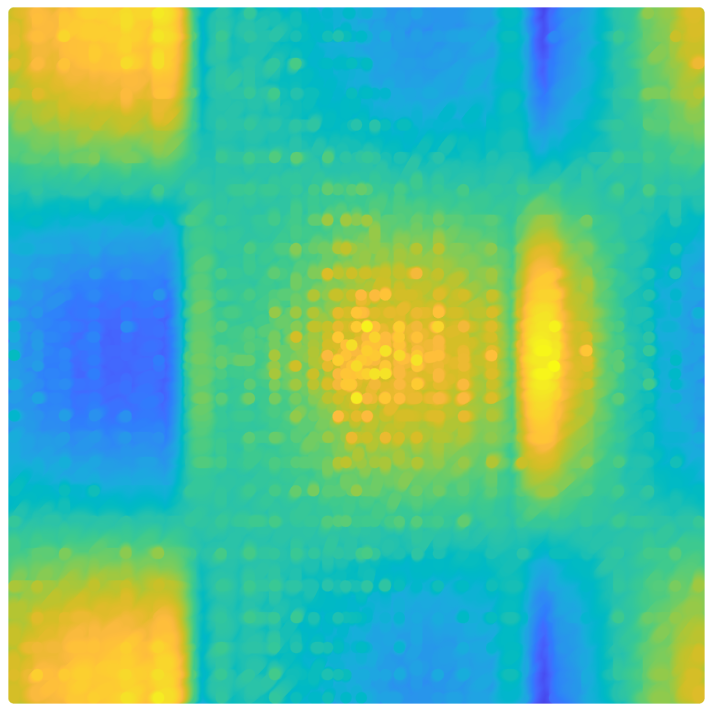}
	}
	\subfloat[WADG (penalty flux)]{
		\includegraphics[width=0.45\linewidth,height=0.3\linewidth]{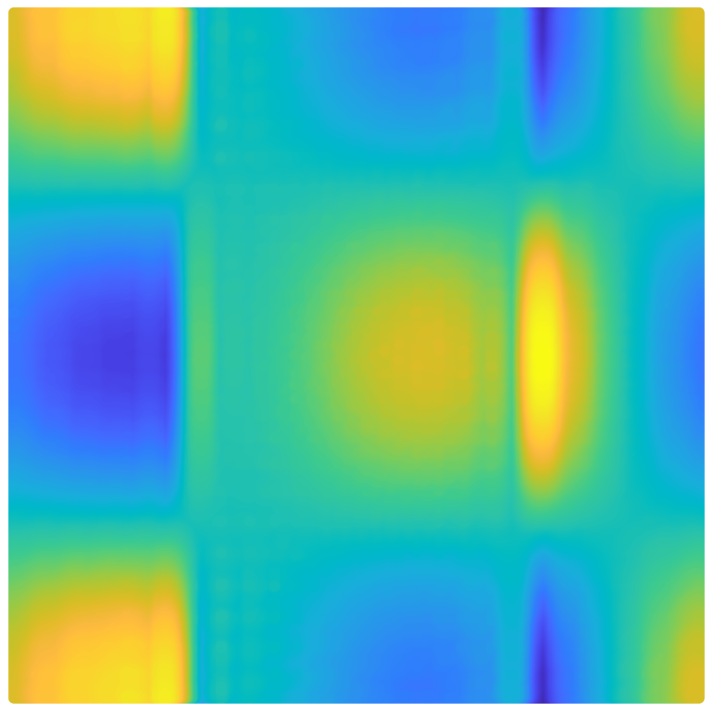}
	}

	\caption[Visuallization of discontinuous solutions using DG and WADG]{Visuallization of discontinuous solutions using DG and WADG}
	\label{fig:dis-condition}
\end{figure}
We first implement with a central flux ($\tau_q=0$) and observe similar convergence rates for WADG methods as previous experiments for discontinuous constant solutions in Figure~\ref{fig:energy-dis-wave}. Then, we use a penalty flux ($\tau_q=1$) and investigate energy dissipation for discontinuous initial conditions. From Figure~\ref{fig:energy-dis-wave-penalty}, we can find that the changes in energy using DG and WADG methods are similar, which is consistent with the observation for sufficiently regular solutions in the previous section. In Figure~\ref{fig:dis-condition}, we present the visualization of discontinuous solutions using DG and WADG methods. We can observe that results given by DG and WADG are almost identical, which suggests that in practice, weight-adjusted DG behaves similarly to DG. Furthermore, the results when using central fluxes for both DG and WADG methods appear to show spurious oscillations in the solution. These oscillations are not present when using penalty fluxes, which suggests that the dissipation introduced by penalty fluxes is sufficient to remove spurious numerical artifacts while retaining accuracy.  

\section{Application to multi-patch DG isogeometric analysis on quadrilateral meshes}\label{sec:ale-spline}
Isogeometric analysis is used to integrate NURBS-based computer representations of geometry directly into the finite element method, with the goal of addressing geometric bottlenecks into the finite element framework \cite{hughes2005isogeometric}. Discretizations using B-spline bases possesses advantages over traditional finite element methods, including the elimination of geometric approximation errors, greater efficiency per degree of freedom, and a larger maximum stable time-step size for explicit time-stepping schemes \cite{chan2018multi}. Here, we apply the proposed ALE-DG formulation on moving quadrilateral meshes using B-spline bases to demonstrate its generality \cite{hughes2008duality,hughes2014finite,evans2009n}. 

Since B-spline bases are non-polynomial, common techniques for moving curved meshes, like mass-lumping, can lose high order accuracy or stability. However, using the weight-adjusted approach, the proposed method is able to deliver high order accuracy and near energy stability for polynomial bases and B-spline bases as well. For the weight-adjusted approximation to the curved mass matrix, high-dimensional B-spline polynomial interpolation and $L^2$ projection operators can be efficiently applied as the Kronecker product of corresponding one-dimensional B-spline operators, while the standard weighted mass matrix in standard DG does not share the same property.

\subsection{B-spline bases}
\begin{figure}[h]
	\centering
	\subfloat[$N=1$]{
		\includegraphics[width=0.32\linewidth]{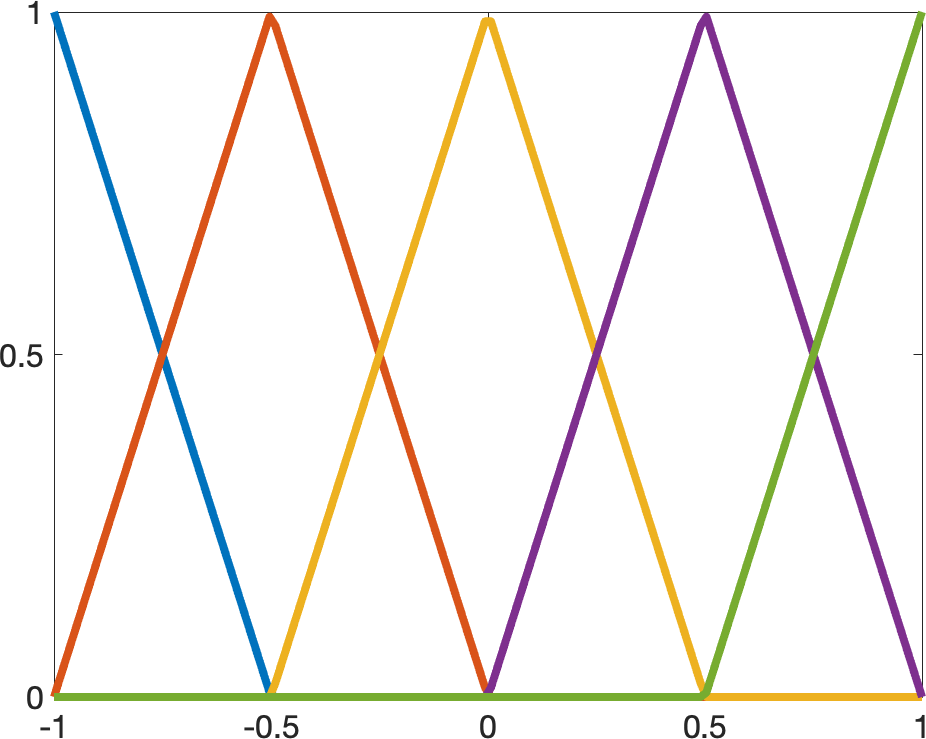}
	}
	\subfloat[$N=2$]{
		\includegraphics[width=0.32\linewidth]{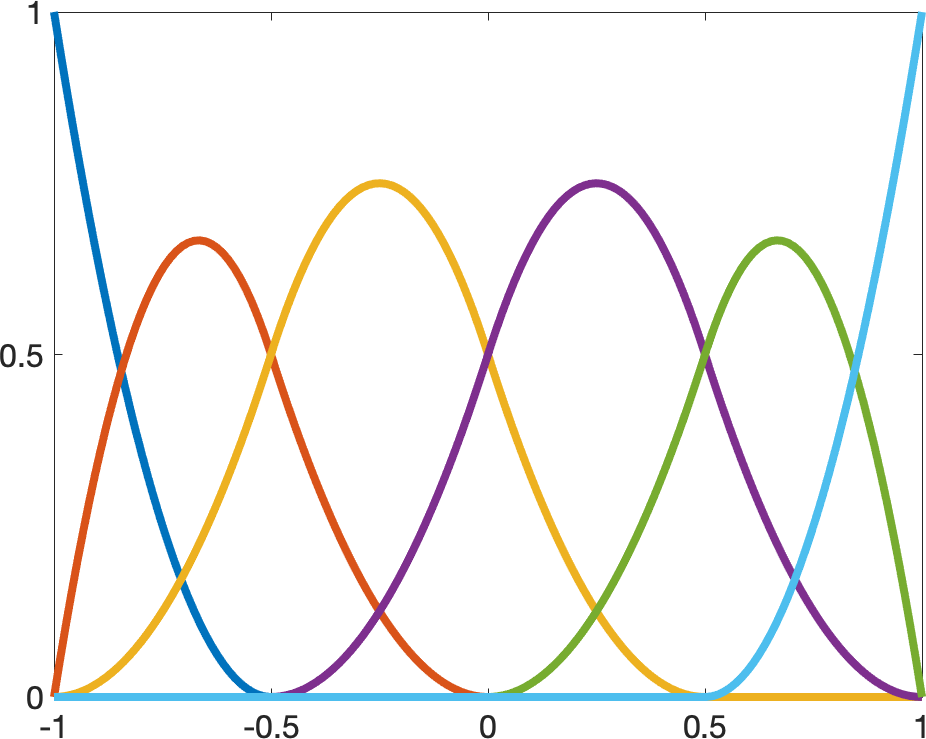}
	}
	\subfloat[$N=3$]{
		\includegraphics[width=0.32\linewidth]{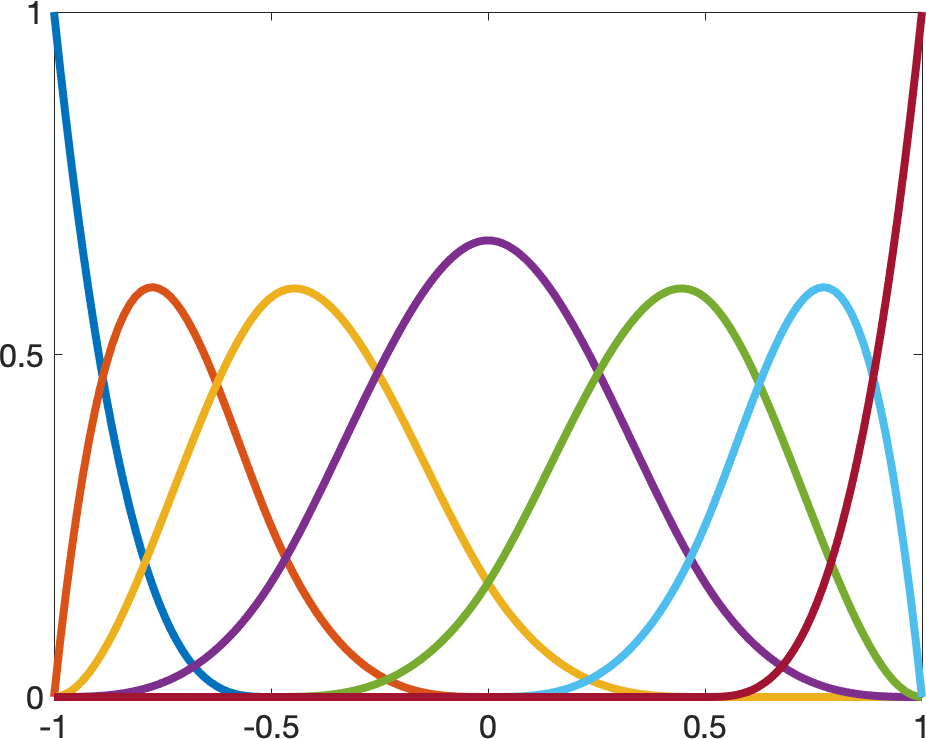}
	}
	\caption[Visuallization of B-spline bases of different degrees.]{Visuallization of B-spline bases of different degrees.}
	\label{fig:bsplins}
\end{figure}

\begin{figure}[h]
	\centering
	\includegraphics[width=0.9\linewidth]{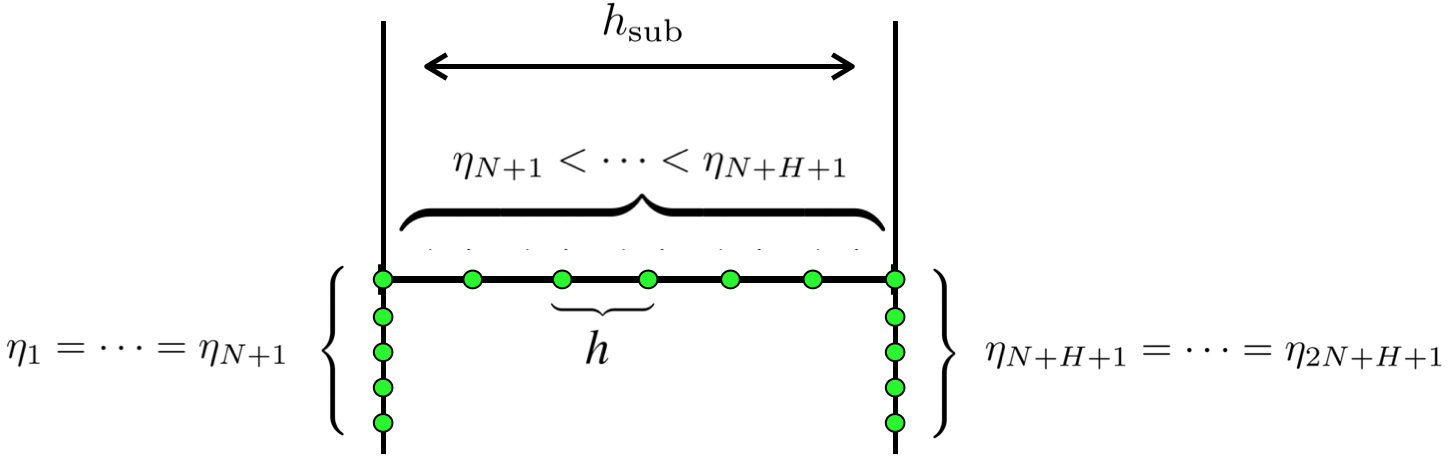}
	\caption{Sketch of the open knot vector on one element.}
	\label{fig:Ksub}
\end{figure}

B-spline basis functions are constructed from an ordered knot vector $\mathbb{H}$
$$\mathbb{H}=\{-1=\eta_1,\dots,\eta_{2N+H+1}=1\},\qquad \eta_i\leq \eta_{i+1}.$$
The unique knots in $\mathbb{H}$ define $H$ non-overlapping sub-elements over which the spline space of degree $N$ is defined. B-spline basis functions (see Figure~\ref{fig:bsplins}) are constructed recursively through the Cox-de Boor recursion formula
$$B_i^0\left(x\right)=\begin{dcases*}
1,&$\eta_i\leq x\leq \eta_{i+1}$,\\
0,&otherwise,
\end{dcases*}\qquad B^k_i\left(x\right)=\frac{x-\eta_i}{\eta_{i+N}-\eta_i}B_i^{k-1}\left(x\right)+\frac{\eta_{i+N+1}-x}{\eta_{i+N+1}-\eta_{i+1}}B_{i+1}^{k-1}\left(x\right).$$
In the case when the denominator vanishes
$$\eta_{i+N}-\eta_i=0\quad \textmd{or}\quad \eta_{i+N+1}-\eta_{i+1}=0,$$
we set the corresponding ratios via $$\frac{x-\eta_i}{\eta_{i+N}-\eta_i}=0\quad\textmd{or}\quad\frac{\eta_{i+N+1}-x}{\eta_{i+N+1}-\eta_{i+1}}=0.$$ The resulting B-spline basis functions have local support and span a piecewise polynomial space of degree $N$. In this work, we restrict ourselves to  the open knot vectors (see Figure~\ref{fig:Ksub}) where the first and last knots are repeated $N$ times
$$\eta_1=\cdots=\eta_{N+1},\qquad \eta_{N+H+1}=\cdots=\eta_{2N+H+1}.$$
We denote the size of sub-elements in B-spline experiments by $h=h_{\text{sub}}/H$,
where $h_{\text{sub}}$ is the maximum distance between spline knots, i.e., the element size.
In case of open knots vectors, the resulting B-spline base contains exactly $N+H$ basis functions. In particular, for $H=1$, this B-spline base recovers the degree $N$ Bernstein polynomial base. We define the one-dimensional degree $N$ B-spline approximation space over the reference interval $\wh{D}=[-1,1]$ as
$$V_h\left(\wh{D}\right)=\text{span}\big\{B^N_i\left(x\right)\big\}^{N+H}_{i=1}.$$

We can extend one-dimensional B-spline bases to $d$ dimensions through a tensor product construction. For example, B-spline basis functions $B^N_{ij}\left(x_1,x_2\right)$ in two dimensions and $B^N_{ijk}\left(x_1,x_2,x_3\right)$ can be defined as
$$B^N_{ij}\left(x_1,x_2\right)=B^N_i\left(x_1\right)B^N_j\left(x_2\right),\quad B^N_{ijk}\left(x_1,x_2,x_3\right)=B^N_i\left(x_1\right)B^N_j\left(x_2\right)B^N_k\left(x_3\right),$$ where $1\leq i,j,k\leq N+H.$ The $d$-dimensional local B-spline approximation space of order $N$ is given by
$$V_h\left(\wh{D}\right)=\text{span}\big\{B^N_i\left(\bm{x}\right)\big\}_{i=1}^{\left(N+H\right)^d},$$
where $B^N_i{\bm{x}}$ denotes the $i$th tensor product spline in $d$ dimensions. 

\subsection{Numerical experiments under B-spline bases}\label{sec:numerical-spline}
In this section, we discuss the extension of the proposed method to quadrilateral meshes using B-spline bases. We start with the energy conservation test for the mesh motion, then compute rates of convergence for wave propagation using analytic solutions. For fair comparisons between results given by B-spline and polynomial bases, we fix the number of ``macro'' elements ($K=2$ for 1D and $K=16$ for 2D) for B-spline experiments and increase the number of sub-elements. The mesh motions for 1D and 2D are given by
\begin{align*}
x&=\xi+\frac{1}{4}\sin\left(\pi \tau\right)\left(1-\xi\right)\left(1+\xi\right),
\end{align*}
and
\begin{align*}
x_1&=\xi_1 + \frac{1}{4}\sin\left(\pi \tau\right)\sin\left(\pi \xi_1\right)\left(1-\xi_1\right)\left(1+\xi_1\right),\\ x_2&=\xi_2 + \frac{1}{4}\sin\left(\pi \tau\right)\sin\left(\pi \xi_2 \right)\left(1-\xi_2\right)\left(1+\xi_2\right),
\end{align*}
separately. We first examine the energy variation for constant solutions on moving meshes for $\p u/\p t=0$. From Figure \ref{fig:spline-1d-mesh-energy} and \ref{fig:spline-2d-mesh-energy}, we observe that the change in energy \textcolor{black}{converges} to zero at a rate of $O(h^{2N+2})$ in both 1D and 2D , as predicted by Theorem \ref{thm:upper-bound}. Furthermore, we investigate energy conservation for wave simulations with central fluxes, and results are given in Figure \ref{fig:energy-central-homo-spline} and \ref{fig:energy-penalty-homo-spline}.  We find that the change in energy converges to zero faster for B-splines than for polynomial bases.
\begin{figure}[h]
	\centering
	\hspace{-4em}
	\subfloat[1D]{
		\begin{tikzpicture}
		\begin{loglogaxis}[
		width=0.55\textwidth,
		height=0.5\textwidth,
		title style = {font=\large},
		xlabel= Mesh size $h$,
		ylabel= $\Delta E$,
		label style = {font=\scriptsize},
		ticklabel style = {font=\scriptsize},
		ymax = 1e-03,
		ymin = 1e-15,
		legend style={font=\tiny},
		legend pos = north west
		]
		
		%greville T=0.5
		\addplot[color=red,mark=o] coordinates {
			(0.5,1.3364e-04)
			(0.25,3.0014e-06)
			(0.125,3.4218e-08)
			(0.0625,4.8457e-10)
		};
		
		\addplot[color=blue,mark=diamond] coordinates {
			(0.5,1.0743e-05)
			(0.25,1.2409e-07 )
			(0.125,2.9302e-10)
			(0.0625,1.0176e-12)
		};
		\addplot[color=red,mark=triangle] coordinates {
			(0.5, 5.9546e-07 )
			(0.25,7.6414e-09)
			(0.125, 3.9582e-12)
			(0.0625,3.7748e-15)
		};

		\logLogSlopeTriangle{0.18}{0.08}{0.47}{6.15}{red};
		\logLogSlopeTriangle{0.18}{0.08}{0.25}{8.19}{blue};
		\logLogSlopeTriangle{0.18}{0.08}{0.05}{10.35}{red};
		\legend{WADG $N=2$,WADG $N=3$,WADG $N=4$}
		\end{loglogaxis}
		\end{tikzpicture}
		\label{fig:spline-1d-mesh-energy}	
	}
	\subfloat[2D]{
		\begin{tikzpicture}
		\begin{loglogaxis}[
		width=0.55\textwidth,
		height=0.5\textwidth,
		title style = {font=\normalsize},
		xlabel= Mesh size $h$,
		label style = {font=\scriptsize},
		ticklabel style = {font=\scriptsize},
		ymax = 1e-03,
		ymin = 1e-15,
		legend style={font=\tiny},
		legend pos = north west
		]
	
		%greville N=2 T=0.5 CFL=0.05
		\addplot[color=red,mark=o] coordinates {
			(0.25,3.1309e-04)
			(0.125,5.9718e-06)
			(0.0625,5.9024e-08)
			(0.03125,8.2192e-10)
			
		};

		%greville N=3 T=0.5 CFL=0.05
		\addplot[color=blue,mark=diamond] coordinates {
			(0.25,1.2553e-05)
			(0.125,3.3387e-07)
			(0.0625,6.4666e-10)
			(0.03125, 2.0799e-12)
			
		};
	
		%greville N=4 T=0.5 CFL=0.05
		\addplot[color=red,mark=triangle] coordinates {
			(0.25, 1.6872e-06)
			(0.125,3.0669e-08 )
			(0.0625,1.0963e-11)
			(0.03125,1.0103e-14)
			
		};

		\logLogSlopeTriangle{0.18}{0.08}{0.49}{6.10}{red};
		\logLogSlopeTriangle{0.18}{0.08}{0.28}{8.30}{blue};
		\logLogSlopeTriangle{0.18}{0.08}{0.085}{10.32}{red};
		\legend{WADG $N=2$,WADG $N=3$,WADG $N=4$}
		\end{loglogaxis}
		\end{tikzpicture}
		\label{fig:spline-2d-mesh-energy}
	}
	\caption[Energy variation of DG and WADG for different orders of approximation (spline)]{Energy variation of DG and WADG for different orders of approximation (spline)}
\end{figure}
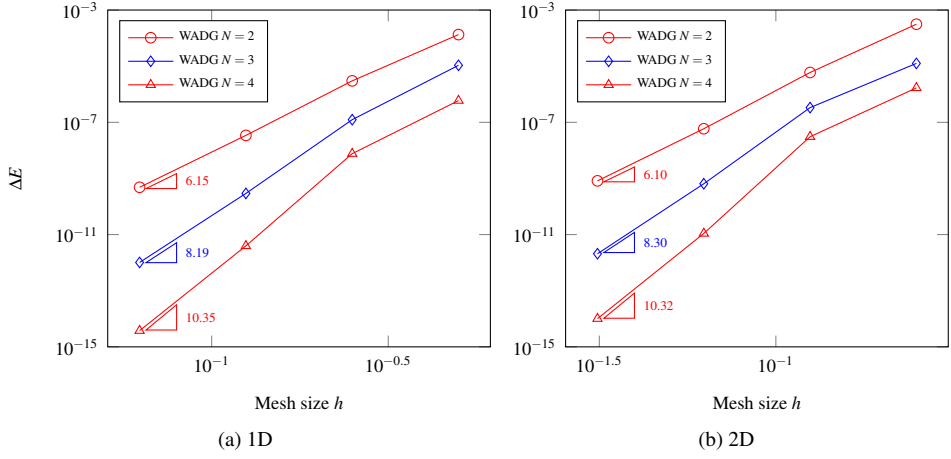

\begin{figure}[h]
	\centering
	
	\hspace{-4em}
	\subfloat[Central flux ($\tau_q=0$)]{
		\begin{tikzpicture}
		\begin{loglogaxis}[
		width=0.55\textwidth,
		height=0.5\textwidth,
		title style = {font=\large},
		xlabel= Mesh size $h$,
		ylabel= $\Delta E$,
		label style = {font=\scriptsize},
		ticklabel style = {font=\scriptsize},
		ymax = 1e-02,
		ymin = 1e-16,
		legend style={font=\tiny},
		legend pos = north west
		]
		
		%greville N=2 T=1.5
		\addplot[color=red,mark=o] coordinates {
			(0.25,2.8742e-04)
			(0.125, 6.7951e-05)
			(0.0625,6.1171e-07)
			(0.03125,8.4653e-09)
		};
	
		%greville N=3 T=1.5
		\addplot[color=blue,mark=diamond] coordinates {
			(0.25,   1.7857e-05)
			(0.125, 2.7574e-07)
			(0.0625,7.7906e-10)
			(0.03125,2.5498e-12)
		};

		% greville N=4 T=1.5
		\addplot[color=red,mark=triangle] coordinates {
			(0.25,  1.2266e-06)
			(0.125, 3.2541e-08)
			(0.0625,1.1542e-11)
			(0.03125,1.0658e-14)
		};

		\logLogSlopeTriangle{0.18}{0.08}{0.56}{6.22}{red};
		\logLogSlopeTriangle{0.18}{0.08}{0.31}{8.30}{blue};
		\logLogSlopeTriangle{0.18}{0.08}{0.145}{10.41}{red};
		\legend{WADG $N=2$,WADG $N=3$,WADG $N=4$}
		\end{loglogaxis}
		\end{tikzpicture}
		\label{fig:energy-central-homo-spline}	
	}
	\subfloat[Penalty flux ($\tau_q=1$)]{
		\begin{tikzpicture}
		\begin{loglogaxis}[
		width=0.55\textwidth,
		height=0.5\textwidth,
		title style = {font=\large},
		xlabel= Mesh size $h$,
		label style = {font=\scriptsize},
		ticklabel style = {font=\scriptsize},
		ymax = 1e-02,
		ymin = 1e-16,
		legend style={font=\tiny},
		legend pos = north west
		]
	
		%greville N=2 T=1.5
		\addplot[color=red,mark=o] coordinates {
			(0.25,  4.2025e-03 )
			(0.125, 2.5236e-05 )
			(0.0625,1.3749e-07  )
			(0.03125,1.2545e-09)
		};

		%greville N=3 T=1.5
		\addplot[color=blue,mark=diamond] coordinates {
			(0.25,  3.0314e-04)
			(0.125,  1.1570e-06)
			(0.0625,7.4772e-10)
			(0.03125,1.4314e-12)
		};
	
		%greville N=4 T=1.5
		\addplot[color=red,mark=triangle] coordinates {
			(0.25,1.9230e-05 )
			(0.125, 2.5481e-07)
			(0.0625,1.6330e-11)
			(0.03125,4.4409e-16)
		};
		
		\logLogSlopeTriangle{0.18}{0.08}{0.5}{6.78}{red};
		\logLogSlopeTriangle{0.18}{0.08}{0.3}{9.03}{blue};
		\logLogSlopeTriangle{0.18}{0.08}{0.055}{15.1}{red};
		\legend{WADG $N=2$,WADG $N=3$,WADG $N=4$}
		\end{loglogaxis}
		\end{tikzpicture}
		\label{fig:energy-penalty-homo-spline}
	}
	\caption[Energy variation of DG and WADG for different orders of approximation (spline)]{Energy variation of DG and WADG for different orders of approximation (spline)}
\end{figure}

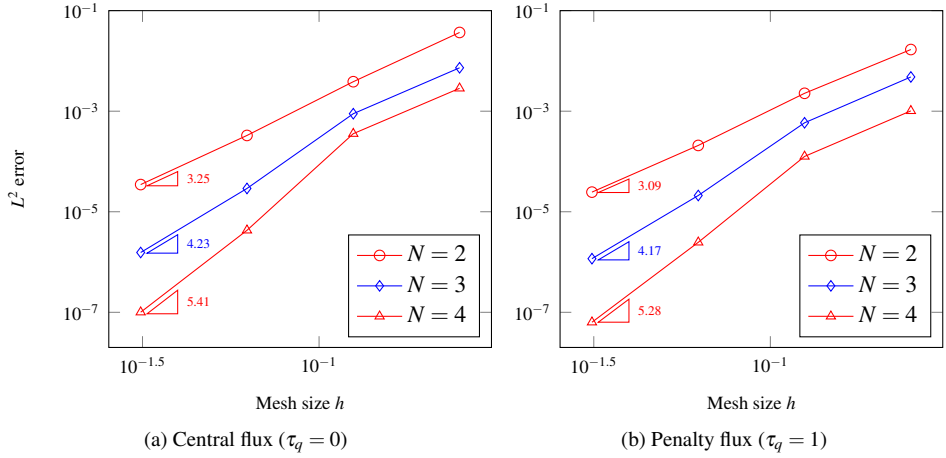
\begin{figure}[H]
	\centering
	\hspace{-4em}
	\subfloat[Central flux ($\tau_q=0$)]{
		\begin{tikzpicture}
		\begin{loglogaxis}[
		width=0.55\textwidth,
		height=0.5\textwidth,
		title style = {font=\large},
		xlabel= Mesh size $h$,
		ylabel= $L^2$ error,
		label style = {font=\scriptsize},
		ticklabel style = {font=\scriptsize},
		ymax=1e-1,
		ymin=0.2e-7,
		legend style={font=\small},
		legend pos = south east
		]
		
		%greville N=2 T=1.5
		\addplot[color=red,mark=o] coordinates {
			(0.25, 3.6810e-02)
			(0.125,3.8602e-03)
			(0.0625,3.2938e-04)
			(0.03125,3.4728e-05)
		};
	
		%greville N=3 T=1.5
		\addplot[color=blue,mark=diamond] coordinates {
			(0.25, 7.2990e-03 )
			(0.125, 8.9118e-04)
			(0.0625, 2.9060e-05)
			(0.03125, 1.5475e-06)
		};

		%greville N=4 T=1.5
		\addplot[color=red,mark=triangle] coordinates {
			(0.25, 2.8526e-03)
			(0.125,3.5715e-04)
			(0.0625,4.2750e-06)
			(0.03125, 1.0045e-07)
		};
		
		\logLogSlopeTriangle{0.18}{0.08}{0.48}{3.25}{red};
		\logLogSlopeTriangle{0.18}{0.08}{0.28}{4.23}{blue};
		\logLogSlopeTriangle{0.18}{0.08}{0.1}{5.41}{red};
		\legend{$N=2$,$N=3$,$N=4$}
		\end{loglogaxis}
		\end{tikzpicture}
		\label{fig:convergence-central-homo-spline}	
	}
	\subfloat[Penalty flux ($\tau_q=1$)]{
		\begin{tikzpicture}
		\begin{loglogaxis}[
		width=0.55\textwidth,
		height=0.5\textwidth,
		title style = {font=\large},
		xlabel= Mesh size $h$,
		%ylabel= $L^2$ error,
		label style = {font=\scriptsize},
		ticklabel style = {font=\scriptsize},
		ymax=1e-1,
		ymin=0.2e-7,
		legend style={font=\small},
		legend pos = south east
		]
	
		%greville N=2 T=1.5
		\addplot[color=red,mark=o] coordinates {
			%(1, 0.505558175415030)
			(0.25,1.6850e-02)
			(0.125, 2.2683e-03)
			(0.0625, 2.0851e-04 )
			(0.03125, 2.4568e-05)
		};

		%greville N=3 T=1.5
		\addplot[color=blue,mark=diamond] coordinates {
			(0.25, 4.7968e-03)
			(0.125, 5.8801e-04 )
			(0.0625,2.1043e-05 )
			(0.03125,1.1662e-06)
		};
		
		%greville N=4 T=1.5
		\addplot[color=red,mark=triangle] coordinates {
			(0.25,1.0153e-03)
			(0.125,1.2593e-04)
			(0.0625,2.4661e-06)
			(0.03125,6.3426e-08)
		};
		
		\logLogSlopeTriangle{0.18}{0.08}{0.46}{3.09}{red};
		\logLogSlopeTriangle{0.18}{0.08}{0.26}{4.17}{blue};
		\logLogSlopeTriangle{0.18}{0.08}{0.075}{5.28}{red};
		\legend{$N=2$,$N=3$,$N=4$}
		\end{loglogaxis}
		\end{tikzpicture}
		\label{fig:convergence-penalty-homo-spline}
	}
	\caption[Convergence behavior for different orders of approximation (spline)]{Convergence behavior for different orders of approximation (spline)}
\end{figure}

Next, we investigate the accuracy of the proposed method by comparing numerical results with analytic solutions of the acoustic wave equation. In Figure \ref{fig:convergence-central-homo-spline} and \ref{fig:convergence-penalty-homo-spline}, we show the rate of convergence to analytic solutions for central fluxes and penalty fluxes. These all converge rapidly with convergence rates close to the theoretical rate $O(h^{N+\frac{1}{2}})$. These observations demonstrate good performance of B-spline bases.

\section{Conclusion and future work}
In this paper, we propose an arbitrary Lagrangian-Eulerian discontinuous Galerkin (ALE-DG) method, which is high order accurate and energy stable up to a term which converges to zero with the same rate as the optimal $L^2$ error estimate. We use the weight-adjusted approach to deal with variable geometries arising from mesh motion. When paired with a skew-symmetric semi-discrete DG formulation, we are able to derive a convergent upper bound for the energy variation. The proposed method is applicable to a wide range of element types and basis functions, which we demonstrate using numerical experiments for both polynomial and B-spline bases on curved triangular and quadrilateral meshes. We numerically verify the energy conservation for the ALE-DG method, and results are consistent with the theoretical upper bound. We also investigate the accuracy of wave simulations and obtain  optimal convergence rates.
\section*{Acknowledgments}
The authors acknowledge the support of the National Science Foundation under awards DMS-1719818 and DMS-1712639.  
%\section*{References}

%%Vancouver style references.
\bibliographystyle{model1-num-names}
\bibliography{refs}

\end{document}